\documentstyle[12pt]{amsart}

\let\Bbb\relax
\newfont{\Bb }{msbm10 scaled 1400} 
\newfont{\Bbb}{msbm10 scaled 1200}  
\newfam\eufam%
\font\euzw=eufm10 scaled 1200%
\font\euac=eufm9%
\textfont\eufam=\euzw\scriptfont\eufam=\euac%
\def\fr{\fam\eufam\euzw}%
\newcommand{\R}{{\Bbb R}}
\newcommand{\C}{{\Bbb C}}
\newtheorem{thm}{Theorem}[section] 
\newtheorem{prop}[thm]{Proposition} 
\newtheorem{cor}[thm]{Corollary} 
\newtheorem{lemma}[thm]{Lemma} 
\newtheorem{dof}{Definition}[section] 
\begin{document}
\title{Homogeneous Special Geometry}
\author{Vicente Cort\'es}
\thanks{Supported by the 
Alexander von Humboldt Foundation, MSRI (Berkeley) and  
SFB 256 (Bonn University).
Research at MSRI is supported in part by grant DMS-9022140.}

\address{\hskip-\parindent Mathematisches Institut der Universit\"at Bonn,
Beringstr. 1, 53115 Bonn, Germany}
\email{vicente@@msri.org {\rm or} vicente@@rhein.iam.uni-bonn.de} 

\curraddr{Mathematical Sciences
Research Institute\\
1000 Centennial Drive\\
Berkeley, CA 94720-5070}

\maketitle 
\pagenumbering{roman}
\begin{abstract} Motivated by the physical concept of special geometry
two mathematical constructions are studied, which relate real hypersurfaces 
to tube domains and complex Lagrangean cones respectively. Me\-thods are 
developed for the classification
of homogeneous Riemannian hypersurfaces and for the classification of 
linear transitive reductive algebraic group actions on pseudo Riemannian
hypersurfaces. The theory is applied to the case of cubic hypersurfaces,
which is the case most relevant to special geometry, obtaining the  
solution
of the two classification problems and the description of the corresponding
homogeneous special K\"ahler manifolds.    
\end{abstract} 
\section*{Introduction}
The concept of special geometry was developed in the context of 
string theory and supergravity. It is a leading principle established
by Witten et al (s.\ e.g.\ \cite{B-W}) that in order to understand and 
construct Lagrangeans $\cal L$ for supersymmetric field theories one has to
interpret the invariance of the action $\int {\cal L} d\mu$ under the gauge
super group as {\it special} geometric property of a Riemannian metric;
namely the metric $g$ defined by the coefficients $g_{ij}$ of the kinetic
terms of the scalar fields $\phi^1,\ldots , \phi^n$ in 
\[ {\cal L} = \sum_{\mu =1}^D\sum _{i,j = 1}^n g_{ij} (\phi^1,\ldots , \phi^n)
\partial_{\mu} \phi^i\partial^{\mu}\phi^j + \cdots \] 

In a more narrow sense (cf.\ e.g.\ \cite{dW-VP0}, \cite{dW-VP} and \cite{St})
{\it special K\"ahler} geometry is the special geometry associated to 
$N = 2$ supergravity-Yang-Mills theories in $D = 4$ spacetime dimensions. 
It has attracted a lot of attention due to its relation to Mirror 
Symmetry, s.\ \cite{Y}.  E.g.\ the Weil-Petersson 
metric of the moduli space of Calabi-Yau 3-folds is special K\"ahler, 
s.\ Example 1 on p.\ \pageref{Ex1}. Similarly, {\it  special real} 
geometry \cite{dW-VP} is defined as the geometry associated to 
$N = 2$ supergravity-Yang-Mills theories in dimension $D = 5$.  
Dimensional reduction of supergravity theories \cite{Cr} from 
dimension $D = 5$ to $D = 4$ induces a correspondence called the 
r-map \cite{dW-VP} relating special real to special K\"ahler
geometry. 

In the first part of the paper, which is a purely mathematical 
presentation of special geometry, we discuss two natural generalizations
of the r-map. 

The basic mathematical object considered is an affine hypersurface
${\cal H} \subset \mbox{\R}^n$ defined by a homogeneous polynomial
$h$, which induces on $\cal H$ a canonical pseudo Riemannian metric,
s.\ \ref{canmetSec}. By the first generalization of the r-map
we canonically associate to $\cal H$ a pseudo K\"ahler tube
domain $U$, s.\ \ref{tubeSec}. The second generalization of the r-map
associates to $\cal H$ a complex Lagrangean cone ${\cal C}\subset T^{\ast}
\mbox{\C}^{n+1}$, s.\ \ref{psKSec}. If $\cal H$ is a {\it cubic} hypersurface,
then the pseudo K\"ahler tube domain $U$ is isomorphic to the projectivized
cone $P({\cal C})$ endowed with a canonical pseudo K\"ahler metric 
induced by the embedding ${\cal C} \subset T^{\ast} \mbox{\C}^{n+1}$, s.\ 
Theorem \ref{gc=gsThm}.  If $\cal H$ is moreover a Riemannian cubic 
hypersurface, then $U \cong P({\cal C})$ is precisely the special 
K\"ahler manifold associated to the special real manifold $\cal H$ 
via the r-map. 

Based on the first construction, we develop a method for the classification
of homogeneous Riemannian hypersurfaces of arbitrary degree $d\ge 2$, 
s.\ \ref{classhomogHSec}. It is easily applied to the cases $2\le d \le 3$.  
As a result, s.\ Theorem \ref{classhomogHThm}, the classification of 
homogeneous cubic Riemannian hypersurfaces and their associated special
K\"ahler (and also quaternionic K\"ahler) manifolds is obtained,
providing a short and conceptual proof for the results of \cite{dW-VP}. 
Our method uses the basic theory of normal J-algebras 
(s.\ Theorem \ref{PSThm}), avoiding the local cubic tensor calculus 
of \cite{dW-VP}. A series of homogeneous Riemannian hypersurfaces 
of arbitrary degree $d\ge 2$ is constructed at the 
end of \ref{classhomogHSec}. 

Finally, a method for the classification of transitive linear reductive
algebraic group actions on pseudo Riemannian hypersurfaces of degree
$d\ge 2$ is presented and carried through for $2\le d \le 3$, s.\ 
\ref{transactSec}. We remark that by Corollary \ref{homogCor} any transitive
linear action of a group $G \subset GL(n,\mbox{\R})$ on a pseudo Riemannian
hypersurface ${\cal H} \subset \mbox{\R}^n$ induces a transitive affine
action of the group $\hat{G} = \mbox{\Bbb R}^n 
\mbox{\Bbb o} (\mbox{\Bbb R}^+ \times G) \subset Aff(\mbox{\Bbb C}^n)$ 
on the corresponding pseudo K\"ahler tube domain $U \subset \mbox{\Bbb C}^n$
by holomorphic isometries. 
If $\cal H$ is cubic, $U$ is special pseudo K\"ahler by Theorem \ref{gc=gsThm}
and we obtain a transitive affine action on the special pseudo
K\"ahler manifold $U$ by holomorphic isometries.    
 
\medskip\noindent 
{\bf Acknowledgements}\\
The author is very grateful to S.-S.\ Chern, R.\ Osserman and MSRI for 
hospitality and support. He would also like to thank D.\ Alekseevsky,
O.\ Baues, J.-M.\ Hwang, M.\ Kontsevich, A.\ Van Proeyen,  E.B.\ Vinberg and 
J.A.\ Wolf for 
discussions related to the subject of this paper.   
\newpage 
\tableofcontents  
\section{Basic constructions of special geometry}\pagenumbering{arabic} 
\subsection{The canonical metric of a hypersurface} \label{canmetSec} 
Let $h$ be a real or complex polynomial in $n$ variables
$x^1,\ldots ,x^n$ which is homogeneous of degree d. 
Consider the level set 
\[ {\cal H}_c(h) = \{ X = (x^1,\ldots ,x^n) \in \mbox{\Bbb K}^n| 
h(X) = c\}\]
where $c\in \mbox{\Bbb K} = \mbox{\Bbb R}$ or {\Bbb C}. If $c\neq 0$, 
we can write ${\cal H}_c(h) = {\cal H}_1(h/c)$. 
\begin{dof} \label{basicDef} A {\bf hypersurface of degree d}  is a 
smooth, open and 
connected subset ${\cal H} \subset  {\cal H}_1(h)$. The 
{\bf basic polynomial} $h$ is homogeneous of degree $d$ and not 
a power of a polynomial of lower degree. 
\end{dof} 
{\bf Remark 1:} Similar constructions to the ones which will be discussed
in the following can be presented for {\em projective} hypersurfaces 
${\cal H} \subset \{ \mbox{\Bbb K}X \in P(\mbox{\Bbb K}^n)| 
h(X) = 0\}$. 

Now we define {\bf the canonical metric} \label{canmetricDef} $g = g(h)$ of a 
hypersurface 
${\cal H} \subset  {\cal H}_1(h)$ of degree $d\ge 2$. 
Let $H\in \vee^d(\mbox{\Bbb K}^n)^{\ast}$ be the symmetric $d$-linear form 
obtained by polarizing $h$, i.e.\ $h(X) = H(X,\ldots , X)$. 
Then we put 
\[ g_{X_0}(X,X) := -(d-1)H(X_0,\ldots , X_0, X,X)\, ,\]
\[ X\in T_{X_0}{\cal H} = \{  X\in \mbox{\Bbb K}^n| H(X_0,\ldots
X_0,X) = 0\}\, .\] 

\noindent 
{\bf  Formula:} 
\begin{equation} \label{HessEqu} 
g = - \frac{1}{d} \partial^2 h =  - \frac{1}{d} \partial^2 \log h 
\quad \mbox{on} \quad {\cal H} \, .\end{equation}  

\noindent
{\bf Remark 2:} By the 
preceding formula,  we can define the canonical metric without assuming 
that $h$ is a homogeneous polynomial. However, we will use the 
homogeneity of $h$ in the next sections, therefore we have  assumed it
from the very beginning.  
\begin{dof} A hypersurface ${\cal H}$ of degree d is called {\bf nondegenerate}
if its ca\-no\-ni\-cal metric is nondegenerate. 
\end{dof} 
In the real case ({\Bbb K} $=$ {\Bbb R}),  the canonical metric $g$ 
defines on   any non\-de\-ge\-ne\-rate hypersurface ${\cal H}$ a canonical  
structure of 
{\bf pseudo Riemannian hypersurface}. 
In the complex case the canonical metric 
defines the structure of  {\bf complex Riemannian hypersurface} on any 
nondegenerate hypersurface.  

\medskip\noindent  
{\bf Remark 3:} {\bf Special real ma\-ni\-folds} in the sense of de Wit 
and Van Proeyen \cite{dW-VP} are precisely Riemannian cubic hypersurfaces.  
   
\subsection{The pseudo K\"ahlerian tube domain associated
to a pseudo Riemannian hypersurface of degree d} \label{tubeSec} 
Let ${\cal H} \subset {\cal H}_1(h) \subset \mbox{\Bbb R}^n$ be a pseudo
Riemannian hypersurface of degree $d$ with canonical metric $g$ 
of signature $(k,l)$. 
We will construct a totally geodesic isometric embedding $\iota : 
({\cal H},g) \hookrightarrow (U, g^c)$ of $\cal H$ into a 
pseudo K\"ahler manifold $(U, g^c)$ of complex dimension
$n$ and complex signature $(k+1,l)$. We consider the positive 
cone ${\cal V} := \mbox{\Bbb R}^+\cdot {\cal H}$  over $\cal H$ and 
the {\bf tube domain} $U = \mbox{\Bbb R}^n + i{\cal V} \subset \mbox{\Bbb C}^n$ 
with complex coordinates $Z = X + i Y$. We define a {\bf canonical} 
pseudo  K\"ahler {\bf metric} $g^c$ on $U$ by the 
\label{potentialP} K\"ahler potential 
\[ K(Z) = -\frac{4}{d} \log (h(Y))\, .\] 
\begin{prop} \label{iotaProp}
 The map $\iota : ({\cal H},g) \hookrightarrow (U, g^c)$
given by $Y \mapsto iY$ is a totally geodesic isometric embedding of the 
pseudo Riemannian hypersurface $({\cal H},g)$ of signature $(k,l)$ 
into the pseudo K\"ahlerian tube domain  $(U, g^c)$ of complex signature
$(k+1,l)$. The cone $i{\cal V} \subset U$ is totally geodesic and isometric 
to the Riemannian product $(\mbox{\Bbb R},can) \times ({\cal H},g)$ via 
\[  \mbox{\Bbb R} \times {\cal H} \ni (t,Y) \mapsto ie^tY \in i{\cal V}\, .\] 
Moreover, $(U, g^c)$ admits the following global isometries preserving the 
cone $i{\cal V}$: 
\begin{itemize}
\item[(i)] scaling $Z\mapsto \lambda Z$ by $\lambda \in  \mbox{\Bbb R}^+$,
\item[(ii)] translations $Z\mapsto Z+X_0$ by real vectors $X_0 \in 
\mbox{\Bbb R}^n$, 
\item[(iii)] reflection $X + iY \mapsto -X + iY$ and 
\item[(iv)] ``inversion''  $X + iY \mapsto X + 
i \frac{Y}{h(Y)^{2/d}}$ with respect to  $i{\cal H}$. 
\end{itemize} 
\end{prop}

\noindent 
{\bf Proof:} By definition we have 
\[ g_Z^c = \sum_{i,j = 1}^n \frac{\partial^2}{\partial z^i \partial \bar{z}^j}
K(Z) dz^i\otimes d\bar{z}^j =\] 
\[ -\frac{1}{d} \sum_{i,j = 1}^n \frac{\partial^2}{\partial y^i \partial 
y^j} \log h(Y) (dx^i \otimes dx^j + dy^i \otimes dy^j)\, , \]  
where $z^j = x^j + iy^j$ are the standard coordinates on $\mbox{\Bbb C}^n 
\supset \mbox{\Bbb R}^n + i{\cal V} = U$, and $Z = X + iY = (x^1,\ldots ,x^n) + i 
(y^1, \ldots , y^n)$. It is clear from this expression, that (ii) and (iii) 
define isometries of $U$ and from the formula (\ref{HessEqu}) it is 
obvious that $\iota : ({\cal H},g) \rightarrow (i{\cal H},g^c|i{\cal H})$ is 
an isometry. Using the formula 
\begin{equation} \label{d2logEqu} \partial^2 \log h = \frac{\partial^2 h}{h} - 
\frac{(\partial h)^2}{h^2}\end{equation}  
and the homogeneity of $h$, it is easy to check that (i) defines an isometric
$\mbox{\Bbb R}^+$ action on $(U, g^c)$ preserving the cone $i{\cal V}$. 
The induced Killing vector field on the cone is just the radial vector field
$R(iY) = iY$. It has unit length and is orthogonal to the hypersurface
$i{\cal H}$. From this it follows that $i{\cal V}$ is a Riemannian product. 
In particular, the inclusions $i{\cal H} \subset i{\cal V} \subset U$ are 
totally geodesic, the cone $i{\cal V}$ being totally geodesic as fixed point 
set of the reflection. Finally, the fact that (iv) is an isometry of $U$
follows from the fact that $iY \mapsto iY/h(Y)^{2/d}$ is an isometry
of $i{\cal V}$. It corresponds to the isometry $(t,Y) \mapsto (-t,Y)$ of 
$\mbox{\Bbb R} \times {\cal H} \cong i{\cal V}$. $\Box$ 

Now we show that our construction behaves nicely with respect to 
group actions. 

\begin{prop} \label{homogpreservProp} Let $G\subset GL(n,\mbox{\Bbb R})$ be a group preserving the
pseudo Riemannian hypersurface ${\cal H} \subset {\cal H}_1(h)
\subset \mbox{\Bbb R}^n$. Then 
$G$ acts on $({\cal H}, g)$ by isometries with respect to the canonical
metric $g$. Moreover, this action on $i{\cal H}\subset U$ is uniquely
extended to a (complex) linear action on $\mbox{\Bbb C}^n \supset U$ 
preserving $U$, namely by the inclusion $G\subset GL(n,\mbox{\Bbb R}) 
\subset GL(n,\mbox{\Bbb C})$. The extended action on $U$ is isometric
(and holomorphic) with respect to the canonical metric $g^c$ on the 
tube domain $U$. 
\end{prop} 

\noindent
{\bf Proof:} Assume that $A \in GL(n,\mbox{\Bbb R})$ preserves 
${\cal H} \subset  \mbox{\Bbb R}^n$. This means that $A^{\ast} h = h$ 
on $\cal H$ and hence on the open set ${\cal V} =  \mbox{\Bbb R}^+\cdot {\cal H}$,
implying $A^{\ast} h = h$. In particular, the transformation $A$ preserves 
the canonical
metric $g$. The action of $A \in GL(n,\mbox{\Bbb R}) 
\subset GL(n,\mbox{\Bbb C})$ on  $\mbox{\Bbb C}^n$ is given by
$X + iY \mapsto AX + iAY$ and it clearly preserves the tube domain $U$ 
and its canonical metric $g^c$. $\Box$ 

From Propositions \ref{iotaProp} and \ref{homogpreservProp} we 
obtain the following corollary. 
\begin{cor} \label{homogCor} If the hypersurface $\cal H$ is homogeneous 
under a subgroup $$G \subset 
 GL(n,\mbox{\Bbb R}) \subset GL(n,\mbox{\Bbb C}),$$ then the 
associated tube domain $(U,g^c)$ is homogeneous as pseudo 
K\"ahler manifold under the subgroup 
\[ \hat{G} = \mbox{\Bbb R}^n 
\mbox{\Bbb o} (\mbox{\Bbb R}^+ \times G) \subset Aff(\mbox{\Bbb C}^n) =
\mbox{\Bbb C}^n \mbox{\Bbb o} GL(n,\mbox{\Bbb C})\, .\] 
\end{cor} 
 
\subsection{Special pseudo K\"ahler geometry and the r-map}
\label{psKSec} 
First we define special pseudo K\"ahler geometry in terms of 
canonical me\-trics on Lagrangean cones. Second we derive the coordinate
expressions (well known in the physical literature) for these 
metrics in terms of a local holomorphic function homogeneous of 
degree 2. Then we define a straightforward generalization of the 
physical r-map. By this we associate a Lagrangean cone $\cal C$ to every real 
hypersurface $\cal H$ of degree d.  If $\cal H$ is a {\em cubic} 
pseudo Riemannian hypersurface, then the projectivized cone $P({\cal C})$ 
with its special metric, s.\ Definition \ref{specialmetricDef}, 
is isometric to the tube domain $(U,g^c)$ associated
to $\cal H$. This shows, in particular, that $(U,g^c)$ is a special 
pseudo K\"ahler manifold. 

Consider the following \label{fadS} {\bf fundamental algebraic data}: 
\begin{itemize} 
\item[1)] A complex symplectic vector space $(V, \omega )$
\item[2)] A compatible real structure, i.e.\ a {\Bbb C}-antilinear 
involution $\tau : V \rightarrow V$ satisfying 
\[ \omega (\tau X, \tau Y) = \overline{\omega (X,Y)}\, ,\quad \forall 
X,Y \in V\, .\] 
\end{itemize} 
Given these data we define a sesquilinear form $\gamma$ on $V$ by
\[ \gamma (X,Y) := i \omega (X, \tau Y) \, .\] 
\begin{lemma} 
The form $\gamma$ is a Hermitian form of signature $(n+1,n+1)$, 
where $\dim_{\mbox{\Bbb C}} V = 2n + 2$. 
\end{lemma} 

\noindent
{\bf Proof:} We first check that $\gamma$ is Hermitian:
\[ \overline{\gamma (X,Y)} = - i\, \overline{\omega (X,\tau Y)} = 
-i\omega (\tau X,Y) = i \omega (Y,\tau X) = \gamma (Y,X)\, .\] 
Now we show that  the Hermitian form $\gamma$ has signature
$(n+1,n+1)$. The restriction of the symplectic form 
$\omega$ to the fixed point set  $V^{\tau}$ of the real structure 
$\tau : V \rightarrow V$ is a real symplectic form on the real vector
space $V^{\tau}$, hence we can choose a basis $p_{\alpha}, q_{\beta}$,
$\alpha ,\beta = 1, \ldots , n+1$, of $V^{\tau}$ such that 
\[ \omega (p_{\alpha}, q_{\beta}) = \delta_{\alpha \beta}\, , \quad 
\omega (p_{\alpha}, p_{\beta}) = \omega (q_{\alpha}, q_{\beta}) = 0\, .\] 
Then  $p_{\alpha}, iq_{\beta}$, $\alpha ,\beta = 1, \ldots , n+1$, is 
a Witt basis for $\gamma$, i.e.\ 
\[ \gamma  (p_{\alpha}, iq_{\beta}) = \delta_{\alpha \beta}\, , \quad 
\gamma (p_{\alpha}, p_{\beta}) = 
\omega (iq_{\alpha}, iq_{\beta}) = 0\, .\]     
This shows that $\gamma$ has signature
$(n+1,n+1)$. $\Box$ 

Remark that up to isomorphism we can assume that the fundamental 
data $(V,\omega , \tau )$ are the following: 
$V = T^{\ast} \mbox{\Bbb C}^{n+1} = T^{\ast} \mbox{\Bbb R}^{n+1} 
+ i T^{\ast} \mbox{\Bbb R}^{n+1}$, $\omega$ the complex bilinear
extension of the standard symplectic form on $T^{\ast} \mbox{\Bbb R}^{n+1}$ 
and $\tau$ complex conjugation with respect to the real form 
$V^{\tau} = T^{\ast} \mbox{\Bbb R}^{n+1}$. 

We recall that a complex  vector subspace $L$ of a complex symplectic 
vector space $(V, \omega )$ is called {\bf Lagrangean} if it is 
maximally isotropic. 

\begin{dof}\label{LagconeDef} Given fundamental algebraic data $(V,\omega , \tau )$ 
a Lagrangean subspace $L \subset V$ is called {\bf nondegenerate} 
if $\gamma |L$ is nondegenerate. A connected submanifold ${\cal C} 
\subset V$ is called  ({\bf nondegenerate}) {\bf Lagrangean cone} 
if $\mbox{\Bbb R}^+ \cdot {\cal C} = {\cal C}$ and if  $L = T_v{\cal C}$ 
is a (nondegenerate) Lagrangean subspace of $V$ for all 
$v\in C$. 
\end{dof}
Remark that a Lagrangean subspace $L \subset V$ is nondegenerate if and only if
$L\cap \tau L = 0$.  
From Definition \ref{LagconeDef} it follows that a nondegenerate Lagrangean 
cone
$\cal C$ is a pseudo K\"ahlerian submanifold of the pseudo K\"ahler 
manifold $(V, \gamma )$. 

\begin{dof} The induced metric 
$g^{{\cal C}} = \gamma |{\cal C}$  is called the {\bf canonical me\-tric} 
of the Lagrangean cone ${\cal C} \subset (V, \omega ,\tau )$. 
\end{dof} 

Now we define a  canonical pseudo K\"ahler metric $g^{P({\cal C})}$ on 
the projective image $P({\cal C}) \subset P(V)$ of a nondegenerate 
Lagrangean cone $\cal C$. For this we assume that $\gamma (u,u) \neq 0$ for 
all $u \in {\cal C}$. If this additional condition is satisfied we shall 
say that the cone $\cal C$ is {\bf properly} nondegenerate. 
Remark that a Riemannian cone $({\cal C},g^{{\cal C}})$ is automatically 
properly nondegenerate. Denote by $\pi : V \rightarrow P(V)$ the 
canonical projection. We define
\[ g^{P({\cal C})}_{\pi u} (d\pi v, d\pi v) = 
\frac{\gamma (v,v)}{\gamma (u,u)} - 
\left|\frac{\gamma (u,v)}{\gamma (u,u)}\right|^2\] 
for $u\in {\cal C} \subset V$, $v\in T_u{\cal C} \subset V$.  

\begin{dof} \label{specialmetricDef} 
The Hermitian metric $g^{P({\cal C})}$ is called the 
{\bf special metric} of the projectivized Lagrangean cone $P({\cal C})$. 
A pseudo K\"ahler manifold is called {\bf special} if it is locally
isometric to $(P({\cal C}),g^{P({\cal C})})$ for some Lagrangean cone 
${\cal C} \subset (V, \omega ,\tau )$. 
\end{dof} 

\begin{prop} \label{univbdlProp} Assume that $\cal C$ is a  properly 
nondegenerate Lagrangean cone.
Then the special metric $g^{P({\cal C})}$ on $P({\cal C})$ is 
a pseudo K\"ahler metric with K\"ahler form $${Im} (g^{P({\cal C})}) = 
2\pi c_1({\cal U},\gamma ),$$ where $c_1 ({\cal U},\gamma )$ is the 
Chern form of the universal bundle $p: {\cal U} \rightarrow P({\cal C})$ 
with Hermitian metric induced by $\gamma$. If the canonical 
metric $g^{{\cal C}}$ has 
Riemannian signature, then the special metric $g^{P({\cal C})}$ is a 
(Riemannian) K\"ahler metric. 
\end{prop} 

\noindent
{\bf Proof:} Let $\zeta$ be a  holomorphic section of ${\cal U}$. Then 
\[ 2\pi c_1({\cal U},\gamma ) = -i \partial \bar{\partial} \log \gamma 
(\zeta , \zeta ) =\] 
\[ -i \frac{\partial \bar{\partial} \gamma (\zeta , \zeta )}{\gamma 
(\zeta , \zeta )} + i\frac{\partial \gamma (\zeta , \zeta ) \wedge 
\bar{\partial} \gamma (\zeta , \zeta )}{\gamma (\zeta , \zeta )^2} =\] 
\[ \frac{{Im}(\gamma )}{\gamma (\zeta , \zeta )} + i\frac{
\gamma(\cdot , \zeta ) \wedge \gamma (\zeta , \cdot )}{\gamma 
(\zeta , \zeta )^2} = {Im}(g^{P({\cal C})}_{p(\zeta )}) \]  
(here we have used the convention $a\wedge \bar{b} = (a \otimes \bar{b} 
- \bar{b} \otimes a)/2$). $\Box$ 

\medskip  
Now we describe Lagrangean cones $\cal C$ and the special metric 
$g^{P({\cal C})}$ on the projectivization $P({\cal C})$ by a 
basic function $F$, thereby relating our presentation to the usual
construction of special K\"ahler geometry in the phy\-si\-cal literature. 
We consider the standard model $V = T^{\ast} \mbox{\Bbb C}^{n+1}$. 
Let $F(Z)$, $Z = (z^0,\ldots , z^n) \in \mbox{\Bbb C}^{n+1}$, be 
a locally defined holomorphic function on $\mbox{\Bbb C}^{n+1}$  
which is homogeneous of degree 2, i.e.\ $F(\lambda Z) = \lambda^2 F(Z)$, 
$\lambda \in \mbox{\Bbb C}-\{ 0\}$,  where this equation makes sense. 
We call $F$ a {\bf basic function}. For example, we may take
$F = p/q$ to be the quotient of two homogeneous polynomials 
$p,q$, $\deg p = \deg q + 2$. We are interested in the image 
${\cal C}_F \subset T^{\ast} \mbox{\Bbb C}^{n+1}$ of the 
differential $dF: Z\mapsto {dF|}_Z$. 

\begin{prop}\label{LagrconeProp} Let $F$ be a basic function. 
Then the connected
components of ${\cal C}_F -\{ 0\}$ are Lagrangean cones. Conversely,
every Lagrangean cone which locally projects isomorphically onto 
$\mbox{\Bbb C}^{n+1}$ is locally of this form for some basic function
$F$.
\end{prop}

\noindent
{\bf Proof:} It is a well known fact in symplectic geometry, that a 
Lagrangean submanifold of  $T^{\ast} \mbox{\Bbb C}^{n+1}$ which 
locally projects isomorphically onto $\mbox{\Bbb C}^{n+1}$ is 
locally the image ${\cal C}_F$ of a  differential $dF$. Now ${\cal C}_F$
is a cone if and only if $dF$ is homogeneous of degree one, i.e. 
if and only if $F$ is homogeneous of degree two. $\Box$ 

Now let ${\cal C} = {\cal C}_F$ be a Lagrangean cone 
which is the image of the differential $dF$ of a basic function $F$. 
We want to express the metrics $g^{\cal C}$ and $g^{P({\cal C})}$
in terms of $F$. Denote by  $Z = (z^0,\ldots , z^n)$ and 
$P = (p_0,\ldots , p_n)$ be the canonical coordinates for 
$T^{\ast} \mbox{\Bbb C}^{n+1}$. 

\begin{prop} \label{canmetricFProp} The canonical metric $g^{\cal C}$ of the Lagrangean cone 
${\cal C} = {\cal C}_F$ at the point ${dF|}_Z\in {\cal C}_F$ is given
by: 
\[ g^{\cal C} = -2 \sum_{j,k = 0}^{n} {Re}\, \left(
i \frac{ \partial^2 F(Z)}{\partial z^j \partial z^k}\right) dz^j \otimes 
d\bar{z}^k = 2 \sum_{j,k = 0}^{n} {Im}\, \left( 
\frac{ \partial^2 F(Z)}{\partial z^j \partial z^k}\right) dz^j \otimes 
d\bar{z}^k\, .\] 
\end{prop} 

\noindent 
{\bf Proof:} The standard symplectic form $\omega$ on 
$T^{\ast} \mbox{\Bbb C}^{n+1}$ is 
\[ \omega = \sum_{j=0}^n (dz^j\otimes dp_j - dp_j\otimes dz^j) \, .\] 
The corresponding Hermitian metric $\gamma = 
i \omega (\cdot , \bar{{\cdot }})$ of signature $(n+1,n+1)$ on  
$T^{\ast} \mbox{\Bbb C}^{n+1}$ is given by 
\[ \gamma = i \sum_{j=0}^n (dz^j\otimes d\bar{p}_j - dp_j\otimes d\bar{z}^j) 
\, .\] 
On ${\cal C} = {\cal C}_F$ we have $p_j = \partial F/\partial z^j$, hence 
\[ g^{\cal C} = \gamma |{\cal C} = i \sum_{j=0}^n (dz^j\otimes d 
\frac{\partial \bar{F}}{\partial \bar{z}^j} - d\frac{\partial F}{\partial z^j}
\otimes d\bar{z}^j) \] 
\[ \stackrel{(\ast )}{=} i \sum_{j,k =0}^n \left( \frac{\partial^2 \bar{F}}{
\partial \bar{z}^j \partial \bar{z}^k}  - 
\frac{\partial^2 F}{\partial  z^j  \partial z^k}\right)dz^j 
\otimes d\bar{z}^k = 2 \sum_{j,k = 0}^{n} {Im}\, \left( 
\frac{ \partial^2 F}{\partial z^j \partial z^k}\right) dz^j \otimes 
d\bar{z}^k\, ,\] 
where at $(\ast )$ we have used that $F$ is holomorphic. $\Box$ 

Since we are assuming that ${\cal C} = {\cal C}_F$, the cone 
${\cal C}\subset T^{\ast} \mbox{\Bbb C}^{n+1}$ is mapped  isomorphically 
onto an open subset of $\mbox{\Bbb C}^{n+1}$  under the canonical 
projection $T^{\ast} \mbox{\Bbb C}^{n+1} 
\rightarrow \mbox{\Bbb C}^{n+1}$. Therefore, we can use  inhomogeneous 
coordinates of $P({\mbox{\Bbb C}}^{n+1})$ as coordinates on $P({\cal C})$, 
e.g.\ $\zeta^j = z^j/z^0$, $j = 1, \ldots , n$. Recall that in order for 
our definition of the special metric $g^{P({\cal C})}$ to make sense  we 
assume that $\cal C$ is properly non degenerate, i.e.\ $\gamma (u,u) \neq 0$ 
for all $u\in {\cal C}$. Without restriction of generality we may assume
$\gamma (u,u) > 0$ for all $u\in {\cal C}$. 

\begin{cor} The special metric $g^{P({\cal C})}$ is given with respect to 
inhomogeneous coordinates on $P^n_{{\mbox{\Bbb C}}}$ by the pseudo 
K\"ahler potential 
\[ K^F(\zeta^1, \ldots ,\zeta^n ) = \log \gamma (\zeta , \zeta ) = 
\log \left( 2 \sum_{j,k = 0}^n {Im}\, \left( \frac{ \partial^2 F(\zeta )}{
\partial z^j \partial z^k}\right) \zeta^j \bar{\zeta}^k \right) \, ,\] 
where $\zeta = ( \zeta^0 = 1, \zeta^1, \ldots , \zeta^n)$. 
\end{cor} 

\noindent 
{\bf Proof:} The first equation follows from Proposition \ref{univbdlProp}, 
since $(\zeta^1, \ldots , \zeta^n) \mapsto {dF|}_{\zeta} \in {\cal C} 
\subset  T^{\ast} \mbox{\Bbb C}^{n+1}$ is a local section of the 
universal bundle over the projectivized cone $P({\cal C})$. Remark that 
$(\zeta^1, \ldots , \zeta^n)$ are the  inhomogeneous coordinates
of the point $\pi {dF|}_{\zeta} \in P({\cal C}) \subset P(T^{\ast} 
\mbox{\Bbb C}^{n+1})$. The second equation follows from Proposition 
\ref{canmetricFProp}. Remark that the functions 
$\frac{\partial^2 F}{\partial z^j \partial z^k}$  are homogeneous of 
degree zero and hence well defined on $P^n_{\mbox{\Bbb C}}$ as they 
must be. $\Box$ 

\medskip 
By Proposition \ref{LagrconeProp} to every local holomorphic  function $F$ 
homogeneous  of degree two we can associate a Lagrangean cone ${\cal C}_F$ 
and the special metric $g^F := g^{P({\cal C}_F)}$ on the projectivization
$P({\cal C}_F)$. The {\bf r-map} maps by definition any homogeneous cubic
polynomial $h(x^1,\ldots x^n)$ on $\mbox{\Bbb R}^n$ to the (possibly 
degenerate) special metric $g^{F_h}$ with 
\[ F_h(Z) = \frac{h(z^1, \ldots , z^n )}{z^0}\, , \quad 
Z = (z^0,\ldots , z^n)\, .\] 
We can easily extend the r-map to homogeneous polynomials $h$ of degree $d$
by 
\[ h\mapsto F_h(Z) =  \frac{h(z^1, \ldots , z^n )}{(z^0)^{d-2}}\, .\]  
Now we want to compare the special metric $g^{F_h}$ on the projectivized 
Lagrangean cone $P({\cal C}_{F_h})$ to the canonical metric $g^c$ on the 
tube domain $U = \mbox{\Bbb R}^n + i {\cal V}$, ${\cal V} = \mbox{\Bbb R}^+ 
\cdot {\cal H}$, associated to a hypersurface ${\cal H} \subset 
{\cal H}_1(h)$. Consider the biholomorphism 
\[ \varphi : \mbox{\Bbb C}^n \supset U \ni Z \mapsto 
\pi {dF_h|}_{(1,Z)} \in P({\cal C}_{F_h}) \subset P(T^{\ast}
\mbox{\Bbb C}^{n+1})\, ,\] 
where $\pi : T^{\ast}\mbox{\Bbb C}^{n+1} \rightarrow P( T^{\ast}
\mbox{\Bbb C}^{n+1})$ is the canonical projection. 

\begin{prop} \label{gsProp} The metric $g^s = \varphi^{\ast} g^{F_h}$ on 
the tube 
domain 
$U$ associated to the homogeneous polynomial $h$ of degree $d$ has the 
K\"ahler potential 
\[ K^s(Z) = (\varphi^{\ast} K^{F_h}) (Z) = - \frac{4}{d} \log 
(2 {Im} \, (-H(Z,\ldots , Z)(d-2) + H(Z,\ldots , Z, \bar{Z})d))\, ,\] 
where $H$ is the polarization of $h$. In particular, 
$g^{F_h}$ is nondegenerate if and only if $\partial \bar{\partial} K^s$ is 
nondegenerate. 
\end{prop}

\noindent 
{\bf Proof:} We have to show that 
\[ (\partial^2 F_h) (\zeta , \bar{\zeta}) = -H(Z,\ldots , Z)(d-2) 
+ H(Z,\ldots , Z, \bar{Z})d\, ,\]
where $\partial^2 F_h  = ( \frac{\partial^2 F_h}{\partial z^i 
\partial z^j})_{i,j = 0,\ldots ,n}$ is the complex Hessian form 
at the point $\zeta = (1,Z)$. We compute: 
\begin{eqnarray*} (\partial^2 F_h) (\zeta , \bar{\zeta}) &=& h(Z) (d-2)(d-1) 
- \partial h \bar{Z}(d-2) \\ 
& & - \partial h Z (d-2) + (\partial^2h) (Z, \bar{Z}) \\ 
 &=&  H(Z, \ldots , Z) (d-2)(d-1) - H(Z,\ldots , Z, \bar{Z})(d-2) d\\ 
 & & -H(Z, \ldots , Z)(d-2)d +  
 H(Z,\ldots , Z, \bar{Z})(d-1)d\\ 
 &=& -H(Z,\ldots , Z)(d-2) + H(Z,\ldots , Z, \bar{Z})d\, .\quad \Box 
\end{eqnarray*} 

\begin{prop} \label{gs=gcProp}The metrics $g^c$ and $g^s$ on the tube domain associated 
to a homogeneous cubic polynomial $h$ coincide. 
\end{prop} 

\noindent
{\bf Proof:} The following lemma shows that the K\"ahler potentials 
defining $g^s$ and $g^c$ coincide up to an additive constant. $\Box$ 

\begin{lemma} Let $H\in \vee^3 (\mbox{\Bbb R}^n)^{\ast}$ be the 
polarization of the homogeneous cubic polynomial $h$ and 
$Z = X + iY \in \mbox{\Bbb C}^n$, then 
\[ {Im}\, ( - H(Z,Z,Z) + 3 H(Z,Z,\bar{Z})) = 4 h(Y) \, .\] 
\end{lemma}

The following theorem is a consequence of Propositions \ref{iotaProp}, 
\ref{gsProp} and \ref{gs=gcProp}. 
\begin{thm} \label{gc=gsThm} The metrics $g^c = g^s =  \varphi^{\ast} g^{F_h}$ and 
$g^{F_h}$  associated to a pseudo Riemannian cubic hypersurface 
$({\cal H},g) \subset {\cal H}_1(h)$ with canonical metric $g$ 
of signature $(k,l)$ are special pseudo K\"ahlerian of complex 
signature $(k+1,l)$. In particular, the r-map maps special real 
manifolds (i.e.\ Riemannian cubic hypersurfaces) to special 
K\"ahler manifolds.
\end{thm}

\medskip\noindent 
{\bf Example 1:} \label{Ex1} Let $X$ be a K\"ahler manifold with holonomy algebra
${\fr hol} = {\fr su}(3)$, i.e.\ a general Calabi-Yau 3-fold. Consider
the complex vector space $V = H^3(X,\mbox{\Bbb C})$ with standard real 
structure $\tau$,  $V^{\tau} = H^3(X,\mbox{\Bbb R})$. The 
``intersection'' form $\omega$   
\[ \omega (\xi , \eta ) = \int_{X} \xi \wedge \eta \, , \quad 
\xi , \eta \in V\, ,\] 
is a complex skew symmetric bilinear form on $V$ compatible with $\tau$, i.e.\
\[ \overline{ \omega (\xi , \eta )} = \omega (\tau \xi , \tau \eta )\, , 
\quad  \xi , \eta \in V\, .\] 
By the Hodge-Riemann bilinear relations for primitive cohomology, s.\ e.g.\ 
\cite{We} Ch.\ V Sec.\ 6, the 
form $\omega$  is non degenerate and thus a complex symplectic form. 
Indeed, the third cohomology is primitive, i.e.\ 
$\Omega \wedge \xi = 0 \in H^5(X,\mbox{\Bbb C})$ for all 
$\xi \in H^3(X,\mbox{\Bbb C})$, where $\Omega$ is the K\"ahler class of $X$. 
This follows from the formula $h_0^{p,q} = h^{p,q} - h^{p-1,q-1}$, 
which relates the primitive Dolbeault numbers $h^{p,q}_0$ to the usual
ones $h^{p,q}$ and from the fact that the only holomorphic 
form on $X$ (up to scaling) is the volume form, due 
to our holonomy assumption. 

Summing up, to every general Calabi-Yau 3-fold we have  associated the  
fundamental algebraic data  $(V,\omega , \tau )$, cf.\   definition on 
p.\ \pageref{fadS}. 

Again by the Hodge-Riemann bilinear relations, the complementary 
subspaces $W = H^{3,0}(X) + H^{1,2}(X)$ and $\overline{W} := \tau W$ are Lagrangean
(with respect to $\omega$), the Hermitian form $\gamma = i \omega (\cdot ,
\tau \cdot )$ is positively defined on $W\times  \overline{W}$ and 
$H^{3,0}(X)$ and $H^{1,2}(X)$ are $\gamma$-orthogonal. 

Now we consider the moduli space of $X$, i.e\ deformations of its complex 
structure.   It is known (s.\ \cite{Ti}, \cite{To}, cf.\ \cite{Kn}) that there exists a 
local universal deformation $\delta : {\cal X} \rightarrow S$, 
$X_s := \delta^{-1}(s)$, $X_{s_0} = X$, of $X$. The base $(S,s_0)$ 
is the Kuranishi moduli space, which can be identified with a 
connected open subset of $H^{2,1}(X)$. 

Since $X$ is K\"ahler the Dolbeault numbers $h^{p,q}(X_s)$ are independent 
of the deformation parameter $s$, by well known results of Kodaira 
and Spencer's deformation theory (combine \cite{We} Ch.\ V Cor.\ 6.6 and 
\cite{M-K} Thm.\ 4.6). After identification of the fibres 
of $\delta : {\cal X} \rightarrow S$ by a diffeomorphism 
${\cal X} \cong S\times X$ we can define the period map
\[ {Per} : S\ni s\mapsto H^{3,0}(X_s) \in  P(H^3(X,\mbox{\Bbb C})) 
= P(V)\, .\] 
The differential $d{Per}_s$ maps $T_sS \cong  H^{2,1}(X)$ isomorphically onto 
$d\pi_p(H^{3,0}(X_s)$ $+$  $H^{2,1}(X_s)) = d\pi_p H^{2,1}(X_s)$, where 
$\pi : V \rightarrow P(V)$ is the canonical projection and 
$0 \neq p\in {Per}(s)= H^{3,0}(X_s)$. Remark that the subspace 
$W_s = H^{3,0}(X_s) + H^{2,1}(X_s)\subset V$ is still Lagrangean, 
because a local deformation of a K\"ahler manifold is again K\"ahler, 
s.\ \cite{M-K} Thm.\ 4.6. This shows that the lines ${Per}(s)
\subset V$, $s\in S$, form a Lagrangean cone ${\cal C} = \cup_s {Per}(s) 
\subset V$. 

The special metric $g^{P({\cal C})}$ on $P({\cal C}) = {Per}(S)$
is negatively defined and hence $-g^{P({\cal C})}$ is  a genuine K\"ahler 
metric. In fact,  $\gamma$ is negatively defined on 
$H^{2,1}(X_s) = {Per}(s)^{\perp} \cap T_p{\cal C}$, $p\in {Per} (s)-\{ 0\}$;
here $\perp$ denotes the $\gamma$-or\-tho\-go\-nal complement. 
The (negatively defined) special ``K\"ahler'' metric  $g^{P({\cal C})}$ is 
known as 
Weil-Petersson metric. 

\noindent 
{\bf Remark 4:} The Weil-Petersson metric  on the Kuranishi moduli
space for general Calabi-Yau $m$-folds of arbitrary dimension $m$ 
can be defined similarly. The reasons why we have concentrated
on the case $m=3$ are the following: first ${Per}(S)$ is always
a totally isotropic cone with respect to the intersection form,
but in general not maximally isotropic in $H^m(X,\mbox{\Bbb C})$, and 
second the intersection form is symmetric if $m$ is even.

\section{Homogeneous case} 
\subsection{Classification of homogeneous  
Riemannian hypersurfaces and corresponding Siegel domains}
\label{classhomogHSec} 
Let ${\cal H} \subset {\cal H}_1(h) \subset {\mbox{\R}}^n$ be a 
pseudo Riemannian hypersurface of degree $d \ge 2$ with basic
polynomial $h$ and canonical metric $g$. 
We define the real algebraic linear group 
\[ {Aut} (h) = \{ \varphi \in GL(n,\mbox{\Bbb R})|  
\varphi^{\ast} h = h \}\, ,\]
which acts naturally on ${\cal H}_1(h)$. 
\begin{dof}\label{homogDef} A pseudo Riemannian hypersurface  ${\cal H} \subset {\cal H}_1(h)$
of degree $d \ge 2$ is said to be {\bf homogeneous} if the connected
component ${Aut}_0 (h)$ acts transitively on $\cal H$. 
\end{dof} 
From Definition \ref{homogDef} and the definition of the canonical metric
$g$ (p.\ \pageref{canmetricDef}) it follows that a homogeneous
pseudo Riemannian hypersurface $({\cal H},g)$ of degree $d \ge 2$ admits
a transitive group of isometries, namely ${Aut}_0 (h)\subset {Isom} 
({\cal H},g)$. By Corollary \ref{homogCor} to any  homogeneous (pseudo)
Riemannian hypersurface  of degree $d$ we have canonically associated a 
homogeneous (pseudo) K\"ahler manifold, which is special (pseudo) K\"ahler if 
$d = 3$. This motivates the interest in the classification problem for
homogeneous such hypersurfaces, which we study in this section. 
   
The first step is to  reduce the classification of homogeneous  
pseudo Riemannian hypersurfaces $({\cal H}, g)$ of degree $d\ge 2$ to the 
case of hypersurfaces admitting a transitive triangular subgroup
${\cal L} \subset {Aut}_0 (h)$. 
By a decomposition theorem for real algebraic groups due to  Vinberg \cite{V} 
we prove now that this reduction is possible if the 
canonical metric $g$  is  Riemannian.

\begin{thm} \label{step1Thm} Let  $({\cal H}, g)$ be a homogeneous 
Riemannian hypersurface 
of degree $d\ge 2$.Then the homogeneous pseudo K\"ahlerian tube domain
$U = {\mbox{\R}}^n + i {\cal V}$ associated to ${\cal H} \subset 
{\cal H}_1(h) \subset {\mbox{\R}}^n$ (cf.\ Corollary 
\ref{homogCor}) is a {\bf Siegel domain of type I}, i.e.\ 
$\cal V$ is convex and does not contain any line. Moreover, $\cal V$ 
is a connected component of ${\mbox{\R}}^n-{\cal H}_0(h)$. Finally, ${Aut}_0 (h)$ 
admits the polar decomposition 
\[  {Aut}_0 (h) = {\cal K}\cdot {\cal L}\, ,\] 
where $\cal K$ is the stabilizer  of a point $v\in {\cal H}$ and is a maximal 
compact connected subgroup and $\cal L$ is a maximal triangular subgroup 
acting simply transitively on $\cal H$. In particular, $\cal H$ is 
contractible. 
\end{thm} 

\noindent
{\bf Proof:} The stated properties of the cone $\cal V$ follow, by well known
arguments of Koszul \cite{K}, from the existence of the  ($\mbox{\R}^+\times
{Aut}_0 (h)$)-invariant closed 1-form $\alpha = - d \log h$ on the 
homogeneous cone $\cal V = \mbox{\R}^+ {\cal H} \subset \mbox{\R}^n$ with 
positively defined Euclidean covariant derivative $\partial \alpha >0$. 
Recall that the 
Hessian bilinear form $\partial \alpha = - \partial^2 \log h$ defines 
up to a positive factor the canonical ($\mbox{\R}^+\times
{Aut}_0 (h)$)-invariant Riemannian  metric $g^{\cal V} \hat{=} g^c|i {\cal V}$ 
of the cone $({\cal V},g^{\cal V})
 \cong (i {\cal V},g^c|i {\cal V}) \subset (U,g^c)$, cf.\ Propositions 
\ref{iotaProp} and \ref{homogpreservProp}. 

The polar decomposition follows from Vinberg's theorem \cite{V} and the 
convexity of the cone. In fact, since ${Aut} (h)$ is defined by a polynomial 
equation, it is a real algebraic linear group. Hence, by \cite{V}  we have 
${Aut}_0 (h) = {\cal K}\cdot {\cal L}$, where $\cal K$ is maximal 
compact and connected and $\cal L$ is  maximal triangular. 

From the convexity of the cone $\cal V$ it follows that $\cal K$ fixes a point
$v\in {\cal H}$ (centre of gravity of a compact orbit), hence 
$\cal K$ is contained in the isotropy group  ${\cal K}^1$ of ${Aut}_0 (h)$ at 
$v$. Since ${Aut}_0 (h)$ acts by isometries of the Riemannian metric $g$,
the isotropy group ${\cal K}^1$ is compact and hence ${\cal K} = {\cal K}^1$. 
$\Box$

The second step is to reformulate our classification problem in terms of 
metric Lie algebras. 

\begin{dof} \label{metricLieDef} A {\bf metric Lie group} $({\cal L},g)$ is a Lie group $\cal L$
together with a left-invariant (pseudo) Riemannian metric $g$. Its 
{\bf metric Lie algebra} $({\fr l}, \langle \cdot ,\cdot \rangle )$ consists of the 
Lie algebra ${\fr l} = Lie\, {\cal L}$ together with the scalar product
$\langle \cdot , \cdot \rangle = g_e$, $e \in {\cal L}$ the identity. 
A {\bf (pseudo) 
K\"ahler Lie group} $({\cal L},g,\tilde{J})$ is a metric Lie group 
$({\cal L},g)$ together with a parallel, left-invariant and orthogonal
complex structure $\tilde{J}$. Its {\bf (pseudo) K\"ahler Lie algebra}
$({\fr l}, \langle \cdot ,\cdot \rangle ,J)$ consists of its metric Lie algebra
$({\fr l}, \langle \cdot ,\cdot \rangle )$ together with the orthogonal
complex structure $J := \tilde{J}_e$ on $\fr l$. 
We say that $({\fr l}, \langle \cdot ,\cdot \rangle )$ (resp.\ 
$({\fr l}, \langle \cdot ,\cdot \rangle ,J)$) is metric (resp.\ 
pseudo K\"ahler) Lie 
algebra {\bf for} the pseudo Riemannian (resp.\ pseudo K\"ahler) manifold $M$ 
if it is the metric (resp.\ pseudo K\"ahler) Lie algebra of a metric 
(resp.\ pseudo K\"ahler) Lie group which as pseudo Riemannian (resp.\ 
as pseudo K\"ahler) manifold is isomorphic to $M$. A K\"ahler Lie algebra
$({\fr l}, \langle \cdot ,\cdot
 \rangle ,J)$ is called a {\bf normal J-algebra} 
if $\fr l$ is splittable solvable and if its K\"ahler form 
$\langle \cdot , J\cdot \rangle$ is the differential of a 1-form 
$\omega$ on $\fr l$, i.e.\ if $\omega ([X,Y]) = - \langle X,JY\rangle$ 
for all $X,Y \in {\fr l}$.    
\end{dof} 

It is known \cite{PS} that any normal J-algebra is K\"ahler Lie algebra
for a bounded homogeneous domain, the K\"ahler metric not necessarily 
being the Bergmann metric. 

\noindent
{\bf Example 2:}  The basic examples of normal J-algebras are the following: 
A {\bf key algebra} ${\fr f} = span\{ G,H\}$ {\bf with root} $\mu > 0$ 
is defined in terms of the orthonormal basis $G = JH, H$ by the formula 
$[H,G] = \mu G$. Given $\fr f$ and a Euclidean vector space $\fr x$ with 
orthogonal complex structure, 
$\fr e = f + x$ carries a canonical Euclidean scalar product $\langle \cdot 
,\cdot \rangle$ and complex structure $J$. The structure of {\bf 
elementary  K\"ah\-le\-rian Lie algebra} 
with key subalgebra $\fr f$ is 
defined on $\fr e$ by the formulas
\[ ad_H|{\fr x} = \frac{\mu}{2} Id,\quad ad_G|{\fr x} = 0 \quad 
\mbox{and} \quad [X,Y] = \mu \langle JX,Y\rangle G \quad \mbox{for} 
\quad X,Y\in {\fr x}\, .\] 
The metric Lie algebra $\fr e$ is determined up to isomorphism (i.e.\ 
orthogonal isomorphism of Lie algebras)  
by $n = \dim_{\mbox{\Bbb C}}\, {\fr x}$ and $\mu$. If we wish to 
specify these parameters, we shall write ${\fr e} = {\fr e} 
(n+1,\mu )$. The elementary  K\"ah\-le\-rian Lie algebra ${\fr e} 
(n+1,\mu )$ is K\"ah\-ler Lie algebra for the complex hyperbolic space
$H_{\mbox{\C}}^{n+1}$ with suitably normalized holomorphic sectional
curvature. Remark that  as Lie algebra ${\fr e} (n+1,\mu )$ is isomorphic
to the Iwasawa Lie algebra of the semisimple Lie algebra ${\fr su}(1,n+1)$. 

Let  $({\cal H}, g)$ be a homogeneous Riemannian hypersurface 
of degree $d\ge 2$ with basic polynomial $h$ and ${\cal B}_0 \subset {Aut}_0 (h)$ a simply transitive
triangular subgroup, which exists by Theorem \ref{step1Thm}. By 
Corollary \ref{homogCor} the Lie group ${\cal B} := \mbox{\R}^+ \times 
{\cal B}_0 \subset GL(n,\mbox{\R})$ (resp.\ ${\cal U}_0 := 
\mbox{\R}^n \mbox{\Bbb o}{\cal B}\subset Aff(\mbox{\Bbb C}^n)$) acts
simply transitively and isometrically (resp.\  and isometrically and 
biholomorphically) on the cone $({\cal V},g^{\cal V}) \cong 
(i {\cal V},g^c|i {\cal V})$ (resp.\ on the Siegel domain 
$(U = \mbox{\R}^n + i{\cal V}, g^c)$). The orbit map 
\[ \varphi : {\cal U}_0 \ni u \mapsto u(ip)\in U\, , \quad p\in {\cal H}\, ,\]
induces diffeomorphisms ${\cal B}_0\cong {\cal B}_0(ip) = i{\cal H} \cong 
{\cal H}$, ${\cal B}\cong {\cal B} (ip) = i{\cal V} \cong {\cal V}$ and 
${\cal U}_0\cong {\cal U}_0(ip) = U$ defining the structure of metric
Lie group on ${\cal B}_0$, ${\cal B}$ and ${\cal U}_0$ by the condition
that $\varphi$ be an isometry. Moreover, the pull back of the complex 
structure of the tube domain via $\varphi$ is a left-invariant, parallel
orthogonal complex structure on ${\cal U}_0$. Hence we have on 
${\cal U}_0$ the structure of K\"ahler Lie group. Its K\"ahler Lie algebra
$({\fr u}_0, \langle \cdot , \cdot \rangle ,J)$ contains the (metric) Lie algebras 
${\fr b}_0 = Lie\, {\cal B}_0$ and ${\fr b} = Lie\, {\cal B}$ as metric 
subalgebras. 

\begin{prop} \label{normalProp} The K\"ahler Lie algebra
$({\fr u}_0, \langle \cdot ,\cdot 
 \rangle ,J)$ associated to a homogeneous Riemannian
hypersurface ${\cal H} \subset {\cal H}_1(h)$ of degree $d\ge2$ is a normal 
J-algebra and admits the orthogonal
(semidirect) decompositions:
\[ {\fr u}_0 = {\fr b} {+\!\!\!\!\!\ni} J{\fr b}\, , \quad {\fr b} = 
\mbox{\R}B_0 \oplus {\fr b}_0\, ,\] 
where $J{\fr b}$ is an Abelian ideal and the vector $B_0$ is in the centre 
of $\fr b$ and satisfies the equations
\[ ad_{B_0} JX = ad_X JB_0 = JX \quad \mbox{for all} \quad 
X\in {\fr b}\, .\]
 The polynomial $h_0 = h \circ d \varphi_e|J{\fr b}$ is 
homogeneous of degree $d$ and invariant under the adjoint 
action of ${\cal B}_0$ on the ideal $J{\fr b}\subset {\fr u}_0$. 
Moreover, the differential $d\varphi_e$ of the orbit map $\varphi$ 
maps the ${\cal B}_0$-orbit  ${\cal H}^0 = {\cal B}_0(JB_0) \subset 
{\cal H}_1(h_0) \subset J{\fr b}$ of the vector $JB_0$ diffeomorphically
onto the hypersurface ${\cal H} \subset {\cal H}_1(h) \subset \mbox{\R}^n
\subset T_{ip}U$, $ip = \varphi (e)$. It is an isometry with respect
to the canonical metrics on ${\cal H}^0$ and $\cal H$. 
\end{prop} 

\noindent
{\bf Proof:} The decompositions of  metric Lie algebras 
are induced by the decompositions 
${\cal U}_0 = {\cal B} \mbox{\Bbb n} \mbox{\R}^n$ and ${\cal B} = 
\mbox{\R}^+ \times {\cal B}_0$ of metric Lie groups. 

The definition of the K\"ahler form $\Omega$ of the tube domain
$U$ by the K\"ahler potential $K(Z) = - \frac{4}{d} \log h(Y)$ on p.\ 
\pageref{potentialP} shows that $\Omega$ is exact in the complex of 
${\cal U}_0$-invariant differential forms on $U$. In fact, $\omega 
= \overline{\partial} \log h(Y)$ is an ${\cal U}_0$-invariant 
1-form on the tube domain $U$ and 
$d \omega = \partial \overline{\partial} \log h(Y)$ is proportional
to $\Omega$. Since, by Theorem \ref{step1Thm}, ${\fr u}_0$ is splittable
solvable, this shows that $({\fr u}_0, \langle \cdot ,\cdot \rangle , J)$ is a 
normal J-algebra. 

The last statements follow from the fact that the representations 
${\cal B}_0 \rightarrow GL(J{\fr b})$ and  ${\cal B}_0\hookrightarrow 
GL(\mbox{\R}^n)$ are isomorphic via 
\[d \varphi_e |J{\fr b} : {\fr b} + J{\fr b} = {\fr u}_0 \supset 
J{\fr b} \stackrel{\sim}{\longrightarrow} \mbox{\R}^n \subset 
T_{ip}U = \mbox{\R}^n + i\mbox{\R}^n\, .\quad \Box \]      

\begin{dof} \label{sim} Let $\fr x$, $\fr y$ and $\fr z$ be 
pseudo Euclidean vector spaces. A bilinear map $\psi :{\fr x}\times {\fr y}
\rightarrow {\fr z}$ is said to be 
{\bf isometric}, if 
\[ \langle \psi (X,Y),\psi (X,Y)\rangle =\langle X,X\rangle \langle Y,Y\rangle \] 
for all $X\in {\fr x}$ and $Y\in {\fr y}$.  The number $k=\dim {\fr x}$ is  
the  {\bf order} of the 
isometric map.    
An isometric map  $\psi :{\fr x}\times {\fr y} \rightarrow {\fr
z}$   is called {special} if 
$\dim {\fr y}=\dim {\fr z}$. 

The {\bf transpose} $\psi^t : {\fr x}\times {\fr z} \rightarrow {\fr y}$ 
of an isometric map $\psi :{\fr x}\times {\fr y}
\rightarrow {\fr z}$ is defined by
\[ \langle \psi^t (X,Z),Y\rangle : = \langle Z, \psi (X,Y)\rangle \, , \quad
X\in {\fr x}\, ,\quad  Y\in {\fr y} \quad \mbox{and} \quad Z\in {\fr z}\, .\]
\end{dof} 
Remark that the transpose $\psi^t$ of an isometric map $\psi$ is 
isometric  if and only if $\psi$ (and hence $\psi^t$) is a special isometric 
map. 
  
The following important structure result can be extracted from \cite{PS}, cf.\ \cite{G-PS-V}. 

\begin{thm}\label{PSThm} (S.G.\ Gindikin, I.I.\ Pyatecki\u{\i}-Shapiro, E.B.\ 
Vinberg) Any normal J-algebra $({\fr u}_0, \langle \cdot ,\cdot \rangle ,J)$ has an 
orthogonal semidirect decomposition 
\begin{equation} {\fr u}_0 = {\fr e}_1 + {\fr e}_2 + \cdots + {\fr e}_l
\label{rankEqu} \end{equation}
into elementary K\"ahlerian subalgebras ${\fr e}_j = {\fr f}_j + {\fr x}_j$
with root $\mu_j$, 
$j = 1,2,\ldots , l$. More precisely, $[{\fr f}_i,{\fr e}_j] = 0$ and
$[ {\fr x}_i,{\fr e}_j] \subset {\fr x}_i$ if $i<j$. 

The normal J-algebra $({\fr u}_0, \langle \cdot ,\cdot
 \rangle ,J)$ is 
K\"ahler Lie algebra for a Siegel domain of type I if and only if we have the 
following orthogonal decompositions 
\begin{equation} {\fr u}_0 =  {\fr b} {+\!\!\!\!\!\ni} J{\fr b}\, , \quad 
{\fr b} = 
\mbox{\R}B_0 \oplus {\fr b}_0\, ,\label{decompEqu} \end{equation}  
where $J{\fr b}$ is an Abelian ideal and $B_0 = \sum_{i=1}^l
\frac{1}{\mu_i}H_i$ is in the centre 
of $\fr b$. Under this assumption, we have orthogonal decompositions
\[ {\fr x}_i = \sum_{j=i+1}^l {\fr x}_{ij}\, , \quad 
 {\fr x}_{ij} =  {\fr x}_{ij}^- +  {\fr x}_{ij}^+\, ,\quad 
 {\fr x}_{ij}^+ = J{\fr x}_{ij}^-\, ,\quad i = 1,\ldots ,l-1\, , \quad 
{\fr x}_l = 0\, ,\]
such that  
\[ {\fr b} = {\fr a} + \sum_{j>i}{\fr x}_{ij}^-\, , \quad \mbox{where}\quad  
{\fr a} = span\{ H_i | i= 1,\ldots ,l\}\] 
and we have the commutator 
relations
\[  [ {\fr f}_k,{\fr x}_{ij}] = 0 \quad \mbox{if} \quad i<j ,\; i<k\quad 
\mbox{and} \quad j\neq k\, ;\] 
\begin{eqnarray*} 
ad_{H_j} | {\fr x}_{ij}^{\pm} &=& \pm \frac{\mu_j}{2} {Id}\, , \\  
ad_{G_j} | {\fr x}_{ij}^+ &=& 0\, ,\\
ad_{G_j} | {\fr x}_{ij}^- &=& -\mu_jJ \quad \mbox{if} \quad i<j\, ;
\end{eqnarray*} 
\[ [{\fr x}_{st},{\fr x}_{ij}] = 0\quad \mbox{if} \quad  s<t,\; i<j,\; i<s 
\quad \mbox{and} \quad s\neq j \neq t\, ;\] 
\begin{eqnarray*}
{[}{\fr x}_{jk}^{\pm},{\fr x}_{ij}^-{]} &\subset & {\fr x}_{ik}^{\pm} \, ,\\
{[}{\fr x}_{jk},{\fr x}_{ij}^+{]} &=& 0 \quad \mbox{if} \quad i<j<k\, ;
\end{eqnarray*}    
\begin{eqnarray*}
{[}{\fr x}_{kj}^{\pm},{\fr x}_{ij}^{\mp}{]} &\subset & {\fr x}_{ik}^+\, ,\\
{[}{\fr x}_{kj}^{\pm},{\fr x}_{ij}^{\pm}{]} &=& 0 \quad \mbox{if} 
\quad i<k<j\, . 
\end{eqnarray*}   
The Lie bracket $[ \cdot ,\cdot ]:{\fr x}_{jk}^{\pm}
\times {\fr x}_{ij}^-\rightarrow
{\fr x}_{ik}^{\pm}$ for $i<j<k$ is given by an isometric map 
$\psi_{ijk}:{\fr x}_{jk}^- \times {\fr x}_{ij}^- \rightarrow
{\fr x}_{ik}^-$ as follows: 
\begin{eqnarray*}
{[}X,Y{]} &=& \frac{1}{\sqrt{2}} \psi_{ijk}(X,Y)\, ,\\
{[}JX,Y{]} &=& J{[}X,Y{]}\, ,\quad X\in {\fr x}_{jk}^-\, ,\; Y\in 
{\fr x}_{ij}^-\, . 
\end{eqnarray*}  
The Lie  bracket $[ \cdot ,\cdot ]:{\fr x}_{kj}^{\pm}\times {\fr x}_{ij}^{\mp} 
\rightarrow {\fr x}_{ik}^+$ for $i<k<j$ is given by 
\[ \langle [X,Y],Z\rangle = -\frac{1}{\sqrt{2}}\langle JY ,
 \psi_{ikj}(X,JZ)\rangle\, ,\]  
\[ [X,Y] = [JX,JY]\, , \quad 
X\in {\fr x}_{kj}^-\, ,\; Y\in {\fr x}_{ij}^+\, ,\; Z\in {\fr x}_{ik}^+\, .\]  
\end{thm} 
The number $l$ of elementary K\"ahlerian subalgebras in the decomposition
(\ref{rankEqu})  is called the 
{\bf rank} of the normal J-algebra ${\fr u}_0$. The normal J-algebras for 
Siegel domains of type I will be called normal J-algebras {\bf of type I}. 

\begin{lemma} \label{hyporbitL} Let  $({\fr u}_0, \langle \cdot ,\cdot
 \rangle ,J)$ be a normal J-algebra for a Siegel domain of type I. Then
the orbit ${\cal H} = {\cal B}_0 (JB_0)$ of the vector $JB_0\in J{\fr b}$
under the adjoint action of ${\cal B}_0\subset {\cal U}_0$ on the 
ideal $J{\fr b} \subset {\fr u}_0 = {\fr b} + J{\fr b}$ is a smooth
hypersurface. 
\end{lemma}

\noindent
{\bf Proof:} It is sufficient to prove that the map
\[ {\fr b}_0 \ni X \mapsto ad_XJB_0 \in J{\fr b}\] 
has maximal rank. This follows from the equation $ad_X JB_0 = JX$ ($X\in 
{\fr b}$), which characterizes the vector $B_0 \in {\fr b}$.  
$\Box$ 

By Proposition \ref{normalProp} and Lemma \ref{hyporbitL} the classification 
of homogeneous
Riemannian hypersurfaces of degree $d\ge 2$ reduces to the following
problem, which can be studied using Theorem \ref{PSThm}. 

\noindent
{\bf Problem 1:} \label{prob1} Classify all normal J-algebras 
$({\fr u}_0, \langle \cdot ,\cdot
 \rangle ,J)$ for Siegel domains of type I (i.e.\ which admit the 
decompositions (\ref{decompEqu})) such that  the hypersurface 
${\cal H} = {\cal B}_0 (JB_0)$ is contained in the level set 
${\cal H}_1(h) \subset J{\fr b}$ of a homogeneous polynomial  $h$ of degree 
$d\ge 2$ on  $J{\fr b}$ and the canonical metric of $\cal H$ is 
Riemannian.

Remark that if $d=3$, then to any solution $({\fr u}_0, \langle \cdot ,\cdot
 \rangle ,J)$  
to Problem 1 we can associate  a homogeneous {\em special} 
K\"ahler manifold, namely the K\"ahlerian Siegel domain of type I 
associated to the homogeneous Riemannian {\em cubic} hypersurface $\cal H$. 

For this reason we will give the complete solution to Problem 1 for 
$2\le d\le 3$. Examples of homogeneous Riemannian  hypersurfaces of 
arbitrary degree
$d\ge 2$ and the corresponding homogeneous K\"ahlerian Siegel domains of type I
will be presented at the end of this section. An interesting class of 
homogeneous pseudo Riemannian cubic hypersurfaces and the corresponding
homogeneous special pseudo K\"ahler and pseudo quaternionic  K\"ahler
manifolds will be discussed in \ref{qKSec}.  

Before studying the general case, we consider the normal J-algebras
of the form ${\fr u}_0 = {\fr f}_1 \oplus {\fr f}_2 \oplus \cdots 
\oplus {\fr f}_l$, i.e.\ all the elementary K\"ahlerian subalgebras
are key algebras. We have the orthogonal decomposition ${\fr u}_0
= {\fr b} + J{\fr b}$ with 
\[ {\fr b} = {\fr a} = {span} \{ H_i| 
i= 1,\ldots , l\} = \mbox{\R} B_0 \stackrel{\perp}{\oplus} {\fr b}_0\, , \quad 
B_0 = \sum_{i=1}^l
\frac{1}{\mu_i}H_i\, .\]    

\begin{lemma} \label{homogfctL} If ${\fr u}_0 = {\fr f}_1 \oplus {\fr f}_2 \oplus \cdots 
\oplus {\fr f}_l$ is an orthogonal direct sum of key algebras, then 
a ${\cal B}_0$-invariant homogeneous function $f$ is defined near the 
point $JB_0\in J{\fr b}$ by 
\[ f(\eta ) = \prod_{j=1}^l a_j^{\prod_{k\neq j}\mu_k^2}\, , \quad
\eta = \sum_{j=1}^l a_jG_j\in J{\fr b} = {span}\{ G_j| j=  1,\ldots
 , l\} \, .\] 
Any such (local) function is proportional to a real power of f.   
\end{lemma} 

\noindent
{\bf Proof:} To prove that $f$ is ${\cal B}_0$-invariant it is sufficient to
check that $f$ is annihilated by the adjoint action of $\mu_1 H_1 - \mu_2 H_2
\in {\fr b}_0$. The last statement follows from the fact that a  
homogeneous ${\cal B}_0$-invariant function defined on a neighborhood of 
$JB_0$ in $J{\fr b}$ is uniquely determined by its degree  
$\in \mbox{\R}$ and its constant value on the hypersurface 
${\cal H} = {\cal B}_0(JB_0)$. $\Box$ 

\begin{prop} \label{homogpolProp} Under the assumption of Lemma 
\ref{homogfctL} there exists a 
non constant ${\cal B}_0$-invariant  homogeneous polynomial $h$ on 
$J{\fr b}$ if and only if 
\[ \left( \frac{\mu_i}{\mu_j}\right)^2 \in \mbox{\Bbb Q}\quad \mbox{for all}
\quad i,j = 1, \ldots , l\, .\] 
Under this condition we can, up to homothety of metric Lie algebras, assume
that $\mu_j^2 = p_j/q_j$ is a reduced fraction for all $j$ and 
$p_1 = q_1 = 1$. Then the unique (up to scaling) non constant  
${\cal B}_0$-invariant  homogeneous polynomial of lowest degree $d$ 
is 
\[ h( \eta ) = \prod_{j=1}^l a_j^{(q_j \prod_{k\neq j} p_k)/N} \, , \quad
N = \gcd \{ q_j\prod_{k\neq j} p_k| j = 1, \ldots ,l\} \, ,\]
where ``$\gcd$'' stands for ``greatest common divisor''. The degree $d$ of 
$h$  is 
\[ d = \frac{1}{N} \sum_{j=1}^l (q_j \prod_{k\neq j} p_k)\, .\] 
\end{prop} 

\noindent
{\bf Proof:} A real power of $f$ is a non constant polynomial if and 
only if there exists a $\lambda \in \mbox{\R}^{\ast}$ such that 
\[ \prod_{k\neq j} \mu_k^2 \in \lambda \mbox{\Bbb Q}\quad \mbox{for all}
\quad  j = 1, \ldots , l\]
or equivalently such that for all $i,j = 1, \ldots , l$ we have
\[ \frac{\mu_i^2}{\mu_j^2} = \frac{\prod_{k\neq j} \mu_k^2}{\prod_{k\neq i} 
\mu_k^2} \in \mbox{\Bbb Q}\, .\] 
By scaling the scalar product of the normal J-algebra  
$({\fr u}_0, \langle \cdot ,\cdot \rangle ,J)$ we can assume that 
$\mu_1 = 1$ and hence $\mu_2^2, \ldots \mu_l^2 \in \mbox{\Bbb Q}$. 
Now the remaining statements are immediate. $\Box$ 

\begin{lemma} \label{keyL} The solutions  
$({\fr u}_0, \langle \cdot ,\cdot \rangle ,J)$
to Problem 1 (p.\ \pageref{prob1}) which admit a direct orthogonal 
decomposition ${\fr u}_0 = {\fr f}_1 + {\fr f}_2 + \cdots 
+ {\fr f}_l$ into key algebras and for which $2\le d\le 3$ are 
up to scaling listed in the following table,
where $h(\eta )$, $\eta =  \sum_{j=1}^l a_jG_j$, is a basic polynomial
for the (flat) homogeneous Riemannian hypersurface ${\cal H} = 
{\cal B}_0 (JB_0) \subset {\cal H}_1(h) \subset J{\fr b}$ of degree $d$. 
\begin{tabular}{|c||c|c|c|} \hline 
$d$ & $l$ & ($\mu_i^2\, , \quad i= 1, \ldots , l$) & $h(\eta )$\\ \hline\hline
2 & 2 & (1,1) & $a_1a_2$\\ \hline 
3 & 2 & (1,2) & $a_1^2 a_2$\\ \hline 
3 & 3 & (1,1,1) & $a_1 a_2 a_3$ \\ \hline
\end{tabular} 
\end{lemma}  

\noindent
{\bf Proof:} 
It is straightforward, using Proposition \ref{homogpolProp}, to determine 
the solutions to Problem1 for which 
${\fr u}_0 = {\fr f}_1 + {\fr f}_2 + \cdots + {\fr f}_l$ and $2\le d\le 3$. 
If we normalize $\mu_1 = 1$, then $d = 2$ implies $l =2$, $\mu_2 = 1$ and
$h( \eta ) = a_1a_2$ and $d = 3$ implies either $l = 3$, $\mu_2 = \mu_3 =1$
and $h( \eta ) = a_1a_2a_3$ or $l = 2$ and $\mu_2 \in 
\{ \sqrt{2}, 1/\sqrt{2}\}$. In the last two cases the corresponding 
polynomials are $h( \eta ) = a_1^2a_2$ and $a_1 a_2^2$ respectively. 
Up to scaling the scalar product of ${\fr u}_0$, it is sufficient to consider
the case $\mu_2 = \sqrt{2}$ and $h( \eta ) = a_1^2a_2$. $\Box$ 

\noindent
{\bf Remark 5:} The Riemannian hypersurfaces defined by the two polynomials
$a_1^2 a_2$ and $a_1 a_2^2$ are isometric via the {\em linear}
transformation $a_1G_1 + a_2 G_2 \mapsto a_2G_1 + a_1G_2$. In particular,
the K\"ahlerian tube domains associated to these hypersurfaces are
isomorphic. 

Now we consider the general case of a normal J-algebra for a Siegel domain
of type I. We use the notations and decompositions introduced above, cf.\  
Theorem \ref{PSThm}; in particular, 
\[ {\fr u}_0 = {\fr e}_1 + {\fr e}_2 + \cdots +{\fr e}_l\, , \quad
 {\fr e}_j = {\fr e}(n_j +1, \mu_j)\, , \quad j=1,\ldots , l\, ,\quad n_l = 0
\, .\] 
Consider the decomposition $J{\fr b} = A_{1,0} + A_{0,1}$, where 
\[ A_{1,0}: = J{\fr a} = {span}\{ G_j| j=  1,\ldots
 , l\}\, , \quad A_{0,1} := \sum_{j>i} {\fr x}_{ij}^+  \, .\]  
It defines an $\fr a$-invariant decomposition for the homogeneous
polynomials of degree $d$ on $J{\fr b}$:
\[ \vee^d (J{\fr b})^{\ast} = \sum_{p+q = d} A^{p,q} :=
\pi \left( (A^{1,0})^{\otimes p} \otimes (A^{0,1})^{\otimes q}\right) \, ,\]
where $A^{1,0} \cong A_{1,0}^{\ast}$ and $A^{0,1} \cong  A_{0,1}^{\ast}$ 
are the subspaces of $(J{\fr b})^{\ast}$ which annihilate 
$A_{0,1}$ and $A_{1,0}$ respectively and   $\pi :\otimes^d (J{\fr b})^{\ast}
\rightarrow \vee^d (J{\fr b})^{\ast}$ denotes the natural projection
from the tensor product to the symmetric tensor product. 
We can decompose any homogeneous polynomial $h \in \vee^d (J{\fr b})^{\ast}$
into its pure components: 
\[ h = \sum_{p+q = d} h^{p,q}\, , \quad h^{p,q}\in  A^{p,q} \, .\] 
Put ${\fr a}_0 = {\fr b}_0 \cap {\fr a} = \{ A\in {\fr a}| \langle A,B_0
\rangle = 0\}$. Remark that if $h$ is a ${\fr b}_0$-invariant homogeneous
polynomial, then $h$ and its pure components $h^{p,q}$ are 
${\fr a}_0$-invariant. Next we will apply these general considerations to the 
case of degree $2\le d\le 3$.   

To simplify the terminology we will use the following definition. 
\begin{dof} We say that a (pseudo) K\"ahler Lie algebra 
$({\fr u}_0, \langle \cdot ,\cdot \rangle ,J)$ for a homogeneous 
(pseudo) K\"ahler manifold $U$, s.\ Definition \ref{metricLieDef}, is 
{\bf of degree d} if, up to scaling the metric, $U$ is isometric to 
the (pseudo) K\"ahlerian tube domain associated to a homogeneous
(pseudo) Riemannian hypersurface of degree $d$. 
\end{dof} 
   
\begin{thm} \label{classhomogHThm} The quadratic normal J-algebras are
(up to scaling) precisely the normal J-algebras of type I (s.\ Theorem 
\ref{PSThm})  of the form 
\[ {\fr u}_0 = {\fr e}_1 + {\fr f}_2\, , \quad \mu_1 = \mu_2 =1\, .\]
The ${\fr b}_0$-invariant homogeneous quadratic polynomial on 
$J{\fr b} = J{\fr a} + {\fr x}_{12}^+$, $$J{\fr a} =  {span}\{ G_1,G_2\},$$
is given by
\[ h (\eta ) = a_1a_2 -\frac{1}{2} \langle X , X\rangle \, , \quad 
\eta = a_1G_1 + a_2 G_2 + X\, , \quad X\in {\fr x}_{12}^+\, .\] 
All cubic normal J-algebras have rank $l = 2$ or $3$. The cubic normal
J-al\-ge\-bras of rank $2$ are (up to scaling) precisely the normal J-algebras
of type I of the 
form 
\[ {\fr u}_0 = {\fr e}_1 + {\fr f}_2\, , \quad \mu_1 = 1\, , \; \mu_2 = 
\frac{1}{\sqrt{2}}\, ,\]
\[h (\eta ) = a_1a_2^2  -\frac{1}{\sqrt{2}}a_2\langle X, X\rangle \, .\] 
The cubic normal J-algebras of rank $3$ are (up to scaling) precisely
the normal J-algebras of type I of the form 
\[ {\fr u}_0 = {\fr u}_0 (\psi ) = {\fr e}_1 + {\fr e}_2 + {\fr f}_3\, ,\quad
\mu_1 = \mu_2  = \mu_3 = 1\, ,\] 
determined by an isometric map $\psi : = \psi_{123} : 
{\fr x}_{23}^-\times {\fr x}_{12}^-\mapsto {\fr x}_{13}^-$, which has 
to be special or of order zero, s.\ Definition \ref{sim}.  
The ${\fr b}_0$-invariant homogeneous cubic polynomial on 
$J{\fr b} = J{\fr a} + {\fr x}_{23}^+ + {\fr x}_{13}^+ + {\fr x}_{12}^+$ 
is 
\[ h (\eta ) = a_1a_2a_3 -\frac{1}{2}\sum_{\alpha =1}^3 a_{\alpha} 
\langle X_{\beta \gamma}, X_{\beta \gamma}\rangle + \frac{1}{\sqrt{2}}
\langle \psi (JX_{23},JX_{12}),JX_{13}\rangle
\, ,\]
where $\eta = \sum_{\alpha =1}^3 (a_{\alpha}G_{\alpha} + X_{\beta \gamma})$, 
$X_{\beta \gamma}\in {\fr x}_{\beta \gamma}$, $\beta < \gamma$ and 
$\{ \alpha , \beta , \gamma \} = \{ 1,2,3\}$.      
\end{thm} 

\begin{cor} The quadratic normal J-algebras are  the normal 
J-algebras for the Hermitian symmetric spaces 
$SO_0(2,2+p)/(SO(2)\times SO(2+p))$. The cubic normal J-algebras are precisely 
the normal J-algebras for the special K\"ahler submanifolds of the Alekseevsky
spaces, s.\ \ref{qKSec}.  
\end{cor}

\noindent
{\bf Proof}(of the theorem){\bf :} Let  
$({\fr u}_0, \langle \cdot ,\cdot \rangle ,J)$ be a quadratic normal J-algebra
and $h$ the basic polynomial of the quadratic Riemannian hypersurface ${\cal 
H} = {\cal B}_0 (JB_0)\subset J{\fr b}$. We have the decomposition
\[ h = h^{2,0} + h^{1,1} + h^{0,2}\] 
into ${\fr a}_0$-invariant homogeneous quadratic polynomials. The polynomial
$h^{2,0}$ is nonzero because otherwise the canonical metric of $\cal H$ 
defined by $h$ would be zero on the subspace $J{\fr a}_0 \subset 
J{\fr b}_0 = T_{JB_0}{\cal H}$, which is impossible for a Riemannian metric. 
Hence $h^{2,0}$ defines a ${\fr a}_0$-invariant nonzero homogeneous quadratic 
polynomial on $J{\fr a}$. The canonical metric defined by $h^{2,0}$ on a
hypersurface in $J{\fr a}$ through the point $JB_0\in J{\fr a}$ is Riemannian.
By Lemma \ref{keyL}  this can only occur if $l =2$, i.e.\ 
${\fr u}_0 = {\fr e}_1 + {\fr f}_2$, $\mu_1 = \mu_2 = 1$ (up to scaling) 
and $h^{2,0}(\eta ) = a_1 a_2$ for $\eta = a_1G_1 + a_2G_2$. 

In the cubic case we have 
\[ h = h^{3,0} + h^{2,1} + h^{1,2} + h^{0,3}\] 
and $h^{3,0}$ is a nonzero ${\fr a}_0$-invariant homogeneous cubic 
polynomial. As above, we conclude in this case that $l = 2$ or $l = 3$. 
If $l = 2$ we have (up to scaling) $\mu_1 = 1$ and either $\mu_2 = \sqrt{2}$ 
and $h^{3,0}(\eta ) = a_1^2a_2$ or $\mu_2 = \frac{1}{\sqrt{2}}$ and 
$h^{3,0}(\eta ) = a_1a_2^2$ for $\eta = a_1G_1 + a_2G_2$. If $l = 3$ 
we have (up to scaling) $\mu_1 = \mu_2 = \mu_3= 1$ and  $h^{3,0}(\eta ) 
= a_1a_2a_3$ for $\eta = a_1G_1 + a_2G_2 + a_3G_3$. 

Using the decomposition of $\vee^3(J{\fr b})^{\ast}$ introduced above we show
next that in the first case the polynomial $h^{3,0}(\eta ) = a_1^2a_2$ 
cannot be extended to a ${\fr b}_0$-invariant homogeneous cubic polynomial 
on $J{\fr b}$ unless ${\fr u}_0 = {\fr e}_1 {\in\!\!\!\!\!+} {\fr f}_2 =
{\fr f}_1 \oplus {\fr f}_2$ is a direct sum of key algebras.   

Assume that $h = h^{3,0} + h^{2,1} + h^{1,2} + h^{0,3}$ is 
such an extension,  
denote by $\pi_{p,q} : \vee^3(J{\fr b})^{\ast}\rightarrow A^{p,q}$, 
$p+q =3$, the natural projection and by $Y\mapsto ad_Y^{\ast}$ the 
representation of ${\fr b}_0$ on $\vee^3(J{\fr b})^{\ast}$ induced
by the adjoint representation of ${\fr b}_0$ on the ideal 
$J{\fr b} \subset {\fr u}_0$. Let $Y\in {\fr x}_{12}^- \subset {\fr b}_0$.
Then the equation $ad^{\ast}_Yh=0$ implies
the equations
\begin{eqnarray*} \pi_{2,1} ad^{\ast}_Yh &=& ad^{\ast}_Yh^{3,0} + 
\pi_{2,1} ad^{\ast}_Yh^{1,2} = 0\, ,\\
\pi_{0,3} ad^{\ast}_Yh &=& \pi_{0,3} ad^{\ast}_Yh^{1,2} = 0
\end{eqnarray*} 
for the polynomial 
\[ h^{1,2}(\eta ) = c_1 a_1 q_1(X,X) + c_2 a_2 q_2(X,X)\, , \quad 
\eta = a_1G_1 + a_2 G_2 + X\, , \] 
where $q_1$ and $q_2$ are quadratic forms on  ${\fr x}_{12}^+\ni X$ 
and $c_1$ and $c_2$ are real constants. 
Since $ad^{\ast}_Yh^{3,0}(\eta ) = -2a_1a_2\langle JY,X\rangle$, the 
first equation is satisfied for all $Y\in {\fr x}_{12}^-\neq 0$ if and only if
$h^{1,2} (\eta ) = -\frac{1}{\sqrt{2}} a_1 \langle X, X\rangle$. However,
this implies 
\[ \pi_{0,3} ad^{\ast}_Yh^{1,2} (\eta ) =  -\frac{1}{\sqrt{2}} 
\langle JY,X\rangle \langle X, X\rangle \, ,\]
so the second equation is not satisfied.  This shows that the case $\mu_2
= \sqrt{2}$ is impossible if ${\fr x}_{12}^- \neq 0$. 

To conclude the proof in the case $l=2$ we have to check that  if 
$\mu_2 = 1$ (resp.\ if $\mu_2 = \frac{1}{\sqrt{2}}$) $h^{2,0}(\eta ) 
= a_1 a_2$ (resp.\ $h^{3,0}(\eta ) = a_1 a_2^2$) is extended by $h$ 
given in the theorem to a homogeneous ${\fr b}_0$-invariant 
quadratic (resp.\ cubic) polynomial on $J{\fr b}$, that the canonical
metric of ${\cal H} = {\cal B}_0(JB_0)\subset {\cal H}_1(h) \subset J{\fr b}$
is Riemannian and, more precisely, that that 
$({\fr u}_0, \langle \cdot ,\cdot \rangle ,J)$
is K\"ahler Lie algebra for the K\"ahlerian Siegel domain associated to 
$\cal H$. This is a special case of Proposition \ref{rk=2deg=dProp}. 

Now we show in the case $l =3$ that if  $h = h^{3,0} + h^{2,1} + h^{1,2} + h^{0,3}$ is a 
${\fr b}_0$-invariant polynomial on $J{\fr b}$ and $ h^{3,0} (\eta ) =
a_1 a_2 a_3$, then $h$ must be the polynomial given above. 

It is easy to check that the only ${\fr a}_0$-invariant element of 
$A^{2,1}$ is zero, the ${\fr a}_0$-invariant elements of $A^{1,2}$ 
are spanned by the three polynomials 
\[ a_{\alpha} \langle X_{\beta \gamma},X_{\beta \gamma}\rangle\, ,\] 
and finally, the ${\fr a}_0$-invariant elements $f$ of $A^{0,3}$ are given by
trilinear functions 
\[ t: {\fr x}_{23} \times {\fr x}_{13} \times {\fr x}_{12} 
\rightarrow \mbox{\R}\, ,
\quad \mbox{i.e.} \quad f(X_{23}+X_{13}+X_{12}) = t(X_{23},X_{13},X_{12})\, ,\]
where
\[ X_{\beta \gamma} \in {\fr x}_{\beta \gamma}
\, , \quad \{ \alpha , \beta , \gamma \} = \{ 1,2,3\} \, , \; \beta < \gamma
\, .\]
This implies 
\[ h^{2,1} = 0 \, ,\quad h^{1,2}(\eta ) = \sum_{\alpha = 1}^3 c_{\alpha} 
a_{\alpha}\langle X_{\beta \gamma},X_{\beta \gamma}\rangle \, , \quad 
h^{0,3}(\eta ) = t(X_{23},X_{13},X_{12})\, ,\]
where $c_{\alpha}$ are real constants, $t$ is a trilinear function as above 
and 
\[ \eta = \sum_{\alpha = 1}^3 (a_{\alpha}G_{\alpha} + X_{\beta \gamma})\, ,
\quad X_{\beta \gamma} \in {\fr x}_{\beta \gamma}
\, , \quad \{ \alpha , \beta , \gamma \} = \{ 1,2,3\} \, , \; \beta < \gamma
\, .\]
Now we let $X\in \sum_{i<j}{\fr x}_{ij}^-$ and consider the equations
\begin{eqnarray*} 0 = \pi_{2,1} ad_X^{\ast} h &=& ad_X^{\ast} h^{3,0} 
+ \pi_{2,1} ad_X^{\ast} h^{1,2}\\
0 = \pi_{1,2} ad_X^{\ast} h &=& \pi_{1,2}ad_X^{\ast} h^{1,2} +
\pi_{1,2}ad_X^{\ast} h^{0,3}\\ 
0 = \pi_{0,3} ad_X^{\ast} h &=& \pi_{0,3} ad_X^{\ast} h^{1,2} + \pi_{0,3} 
ad_X^{\ast} h^{0,3}\, .
\end{eqnarray*}
The first equation is satisfied if and only if $c_1 = c_2 = c_3 
= -\frac{1}{2}$. Then the second equation is equivalent to 
\[ t(X_{23},X_{13},X_{12}) = \frac{1}{\sqrt{2}} \langle \psi (JX_{23},JX_{12})
,JX_{13}\rangle\, .\] 
Now the third equation is satisfied only if ${\fr x}_{23}^- = 0$ or 
if $\psi$ is a special isometric map. In fact, if e.g.\ $X \in {\fr x}_{23}^-$
then 
\[ [ X,\sum_{\alpha = 1}^3 (a_{\alpha}G_{\alpha} + X_{\beta \gamma})] = 
a_3 JX + \langle JX,X_{23}\rangle G_2  + \frac{1}{\sqrt{2}}J\psi^t(X,JX_{13})
\] 
and the third equation reads 
\begin{eqnarray*} 0 &=& \frac{1}{2} \langle JX ,X_{23}\rangle \langle X_{13},
X_{13}\rangle + \frac{1}{2}\langle \psi (JX_{23}, \psi^t(X,JX_{13}),JX_{13}
\rangle \\  
 &=&  \frac{1}{2} \langle JX ,X_{23}\rangle \langle X_{13},
X_{13}\rangle + \frac{1}{2}\langle \psi^t(X,JX_{13}),\psi^t(JX_{23},JX_{13})
\rangle \, .
\end{eqnarray*}
If ${\fr x}_{23}^- \neq 0$ we can choose $X = JX_{23} \neq 0$ and the last 
equation shows that  $\psi^t(X,\cdot) : 
{\fr x}_{13}^- \rightarrow {\fr x}_{12}^-$ is injective, hence 
$\psi$ is a special isometric map. 

It only remains to check that the scalar product of the normal J-algebra
${\fr u}_0$ is (up to scaling) induced by the canonical metric of the tube 
domain associated to the hypersurface
${\cal H} = {\cal B}_0 (JB_0) \subset {\cal H}_1(h) \subset J{\fr b}$. 
This is a straightforward computation, cf.\ Proposition \ref{rk=2deg=dProp}. 
Remark that the pure component $h^{0,3}$ involving the special isometric 
map $\psi$ plays no role in this calculation,
since its  Hessian vanishes at the point $JB_0$. $\Box$ 

\medskip
Now we give examples of normal J-algebras of arbitrary degree 
$d = 2, 3,4,\ldots$. For every degree $d\ge 2$ we construct a series  
$({\fr u}_0(p, d-1),\langle \cdot ,\cdot \rangle ,J)$, $p\in \mbox{\Bbb N}_0$,
of normal J-algebras of  degree $d$ and rank $2$. The  subalgebra 
${\fr b}_0$ for these series
is up to scaling isomorphic to the Iwasawa algebra   of ${\fr so}(1,p+1)$ 
with scalar product
induced by the Riemannian metric of hyperbolic 
$(p+1)$-space $H^{p+1}_{\mbox{\Bbb R}} = SO_0(1,p+1)/SO(p+1)$. 

From now on we write  ${\fr u}_0$ instead of 
$({\fr u}_0, \langle \cdot , \cdot \rangle ,J)$. The scalar product  
and complex structure are understood and always denoted by 
$\langle \cdot , \cdot \rangle$ and $J$ respectively.

For every $s=1,2,\ldots $ consider the K\"ahlerian 
Lie algebra 
\[ {\fr u}_0 = {\fr u}_0(p,s) = {\fr e}(p+1,1) {\in\!\!\!\!\!+} 
{\fr e}(1,\frac{1}{\sqrt{s}}) = ({\fr f}_1 + {\fr x}_1) + {\fr f}_2\, ,\] 
where the semidirect orthogonal sum of the key algebra ${\fr f}_2$ with root 
$\mu = \frac{1}{\sqrt{s}}$ and the ideal ${\fr f}_1 + {\fr x}_1$ is 
defined by the condition that 
$ad_{{\fr f}_2}|{\fr x}_1$ has weight decomposition ${\fr x}_1 
= {\fr x}_{12}^- + {\fr x}_{12}^+$, $p = \dim {\fr x}_{12}^-$.

\begin{prop} \label{rk=2deg=dProp}${\fr u}_0 = {\fr u}_0(p,s) 
= {\fr b} + J{\fr b}$, ${\fr b} = 
\mbox{\R}B_0 \oplus {\fr b}_0$, s.\ Theorem \ref{PSThm},   is a 
normal J-algebra
of degree $d = s+1$ with ${\fr b}_0$-invariant polynomial
\[ h(\eta ) = a_1 (\mu a_2)^s - \frac{1}{2}(\mu a_2)^{s-1} 
{\langle X,X\rangle}\, ,\quad \mu = \frac{1}{\sqrt{s}}\, , \]
where $\eta = a_1 G_1 + a_2 G_2 + X$, $X\in {\fr x}_{12}^+$. 
Its  subalgebra ${\fr b}_0$ is 
isomorphic to the Iwasawa algebra of ${\fr so}(1,p+1)$. In particular, 
the Iwasawa subgroup of $SO_0(1,p+1)$ 
acts  simply transitively on a Riemannian hypersurface  
$\cal H$ of degree $d$ and of constant negative curvature. 
The K\"ahlerian Siegel  domain $U$ associated to the Riemannian 
hypersurface ${\cal H} = {\cal B}_0 (JB_0) \subset J{\fr b}$ 
of degree $d$ is symmetric only if $\cal H$ is quadratic and in this case
$U$ is isometric to the Hermitian symmetric space 
$SO_0(2,p+2)/(SO(2)\times SO(p+2))$. If $\cal H$ is cubic (cf.\ Theorem
\ref{classhomogHThm}), then $U$ is special K\"ahler; in fact 
, it is the special 
K\"ahler submanifold of the Alekseevvsky space 
${\cal T}(p)$, s.\ \ref{qKSec}.  
\end{prop} 

\noindent
{\bf Proof:} First we observe that
${\fr u}_0  =  {\fr b} + J{\fr b}$,
where ${\fr b} = span\{ H_1,H_2\} + {\fr x}_{12}^-$ has centre 
$B_0 = H_1 + \frac{1}{\mu } H_2$ if $p\neq 0$. Remark that the metric Lie 
algebra
${\fr b}_0 \subset  {\fr b} = \mbox{\Bbb R} B_0 \stackrel{\perp}{\oplus} 
{\fr b}_0$ is  isomorphic to the Iwasawa algebra
of ${\fr so}(1,p+1)$ with scalar product induced (up to scaling) by the 
hyperbolic metric of $H^{p+1}_{\mbox{\Bbb R}} = SO_0(1,p+1)/SO(p+1)$. 

Consider now the faithful linear representation
$\rho = {ad|} : {\fr b} \rightarrow {\fr gl}(J{\fr b})$. We have 
$\rho (B_0) = {Id}$ and $\rho (B) JB_0 = JB$ for all $B\in {\fr b}$. 
The corresponding representation $R: {\cal B} 
\rightarrow GL(J{\fr b})$ of the simply connected Lie group 
${\cal B} = \mbox{\R}^+ \times {\cal B}_0$ with Lie algebra ${\fr b} = 
\mbox{\R}B_0 + {\fr b}_0$ has the open orbit $R({\cal B})JB_0$ 
and ${\cal H} = {\cal B}_0(JB_0) = R({\cal B}_0)JB_0$ is a codimension one orbit of the 
simply connected Lie group ${\cal B}_0$ with Lie algebra ${\fr b}_0$. 

We will prove that $\cal H$ is defined by a homogeneous polynomial
$h$ of degree $d$ and that the simply connected metric Lie group 
${\cal U}_0$ with metric Lie algebra ${\fr u}_0$  is up to 
scaling isomorphic to the 
K\"ahlerian Siegel domain $U$ associated to the hypersurface $\cal H$. 

For $\eta = a_1G_1 + a_2G_2 + X \in J{\fr b}$, $X \in {\fr x}_{12}^+$, 
we define 
\[ h(\eta ) = a_1 (\mu a_2)^s - \frac{1}{2}(\mu a_2)^{s-1} 
{\langle X,X\rangle}\, ,\]
where $\langle  \cdot , \cdot \rangle$ denotes the Euclidean scalar product of the 
metric Lie algebra ${\fr u}_0$. It is an easy computation to check that $h$ is 
invariant under the linear action of ${\fr b}_0$ on $\vee^d(J{\fr b})^{\ast}$
induced by $\rho |{\fr b}_0 : {\fr b}_0 \rightarrow {\fr gl}(J{\fr b})$. 
Hence ${\cal H} = R({\cal B}_0) JB_0$ is contained in ${\cal H}_1(h)$. 
Let us now consider the orbit map  
\[ \varphi : {\cal B}_0 
\rightarrow {\cal H}\, ,\quad b\mapsto R(b)JB_0 \, .\] 
We check that the Euclidean scalar product $\langle \cdot , \cdot \rangle$ 
on ${\fr b}_0 = T_e{\cal B}_0$ is 
\begin{equation} \label{spE} {\langle \cdot , \cdot \rangle}\quad  = \quad 
\frac{1}{d}(\varphi^{\ast}g)_e \, ,\end{equation}  
where $g = - \frac{1}{d} \partial^2 h$ is the canonical metric 
of the hypersurface $\cal H$ defined by the basic polynomial $h$. We compute  
\[ (\varphi^{\ast}g)_e (H_1 - \mu H_2 ,H_1 - \mu H_2) = 
g_{JB_0}(G_1 -\mu G_2, G_1 -\mu G_2) \]
\[ = g(G_1,G_1) -2\mu g(G_1,G_2) + \mu^2 g(G_2,G_2) \]
\[ = 0 + \frac{1}{d} (2\mu s(\mu \frac{1}{\mu})^{s-1}\mu -
\mu^2s(s-1)  (\mu \frac{1}{\mu})^{s-2}\mu^2)\]
\[ = \frac{1}{d} (2\mu^2 s - \mu^4 s(s-1)) = \frac{1}{d} (2 - 
\frac{s-1}{s}) \]
\[ = \frac{1}{d} (1 + \mu^2) 
= \frac{1}{d}{\langle H_1 - \mu H_2 ,H_1 - \mu H_2\rangle}\quad \mbox{and}\]   
\[ (\varphi^{\ast}g)_e (X_-,X_-) = g_{JB_0}(JX_-,JX_-)\]  
\[ =\frac{1}{d}(\mu \frac{1}{\mu})^{s-1}{\langle JX_-,JX_-\rangle} =  
\frac{1}{d}{\langle X_-,X_-\rangle}\, ,\] 
where $X_- \in {\fr x}_{12}^-$. This proves equation (\ref{spE}). 

Now we extend the linear representation $\rho : {\fr b} \rightarrow 
{\fr gl}(J{\fr b})\subset {\fr gl}( (J{\fr b})\otimes \mbox{\Bbb C})$ 
to an affine representation $\rho : {\fr u}_0 \rightarrow 
{\fr aff}((J{\fr b})\otimes \mbox{\Bbb C})$ of the real Lie algebra
${\fr u}_0$ on the complex vector space $(J{\fr b})\otimes \mbox{\Bbb C} 
= J{\fr b} + i J{\fr b}$. For $B\in {\fr b}$ we define 
\[ \rho (JB) = JB \in  J{\fr b} \subset (J{\fr b})\otimes \mbox{\Bbb C} 
\subset (J{\fr b})\otimes \mbox{\Bbb C} + 
{\fr gl}((J{\fr b})\otimes \mbox{\Bbb C})= 
{\fr aff}((J{\fr b})\otimes \mbox{\Bbb C})\, .\] 
The corresponding extension of $R: {\cal B} 
\rightarrow GL(J{\fr b})\subset GL((J{\fr b})\otimes \mbox{\Bbb C})$ 
will be also denoted by the same letter: 
\[ R: {\cal U}_0 \rightarrow Aff((J{\fr b})\otimes \mbox{\Bbb C})\, .\] 

\noindent
The tube domain $U = J{\fr b} + i {\cal V}$, ${\cal V} = 
\mbox{\Bbb R}^+ {\cal H}\subset J{\fr b}$, is precisely   the orbit of the 
point $iJB_0\in i{\cal V}\subset U$  under the affine representation $R$
of the group ${\cal U}_0$ on $(J{\fr b})\otimes \mbox{\Bbb C}$, that is:  
$U = R({\cal U}_0)(iJB_0)$. 

Now we prove that the scalar product $\langle \cdot , \cdot \rangle$ of the 
metric Lie algebra ${\fr u}_0 = T_e{\cal U}_0$ is 
\[ {\langle \cdot , \cdot \rangle} \quad = \quad \frac{1}{d} (\phi^{\ast} g^c)_e\, ,\]  
where $g^c$ is the canonical metric of the tube domain $U$ 
and 
\[ \phi : {\cal U}_0 \rightarrow U\, ,\quad u\mapsto R(u)(iJB_0)\] 
is the orbit map. 

First we remark that that 
\[ (d\phi )_e : {\fr u}_0 = {\fr b} + J{\fr b} 
\rightarrow T_{iJB_0}U = J{\fr b} + i J{\fr b}\] 
maps the orthogonal subspaces  $\fr b$ and $J{\fr b}$ of ${\fr u}_0$ 
onto the orthogonal subspaces  $iJ{\fr b}$ and $J{\fr b}$ of 
$T_{iJB_0}U$ respectively: 
\[ (d\phi )_e|{\fr b} : {\fr b} \rightarrow iJ{\fr b}\, , \quad
B\mapsto iJB\, \]
\[ (d\phi )_e|J{\fr b} = {Id}: J{\fr b}\rightarrow J{\fr b}\, .\] 
Moreover, $J$ and multiplication by $-i$ are orthogonal endomorphisms 
of ${\fr u}_0$ and $T_{iJB_0}U$ respectively which make the following 
diagram commutative: 
\[ \begin{array}{r@{\,}c@{\,}l@{\:}c@{\;}c@{\,}l} 
    &  {\fr u}_0      &                                  &
\stackrel{(d\phi )_e}{\longrightarrow} &   T_{iJB_0}U        & \\ 
  & \downarrow &  J &                                    & \downarrow &  
-i\\
&   {\fr u}_0     &                                 &
\stackrel{(d\phi )_e}{\longrightarrow}   & T_{iJB_0}U           & 
\end{array} \]
Now since $J$ interchanges the subspaces $\fr b$ and $J{\fr b}$ of 
${\fr u}_0$ and $-i$ interchanges the subspaces $iJ{\fr b}$ and 
$J{\fr b}$ of $T_{iJB_0}U$, it is sufficient to check that the map
$ (d\phi )_e|{\fr b} : {\fr b} \rightarrow iJ{\fr b} \subset T_{iJB_0}U$ 
becomes a linear isometry after scaling it by $\sqrt{d}$. 
We have checked this already on ${\fr b}_0\subset {\fr b} = 
\mbox{\Bbb R} B_0 + {\fr b}_0$, s.\ (\ref{spE}). The vector 
$(d\phi)_eB_0 = iJB_0$ is precisely the radial vector at the point 
$iJB_0\in i{\cal V}$, so it has unit length with respect to the 
canonical metric $g^c$, thanks to Proposition \ref{iotaProp}. 
On the other hand, ${\langle B_0,B_0\rangle} = 1 + \frac{1}{\mu^2} = 1 + s = d$. 
This proves that $(\phi^{\ast}g^c)_e = \frac{1}{d} {\langle \cdot , \cdot 
\rangle}$. 
All remaining statements are easily checked. $\Box$ 

\subsection{Classification of transitive reductive
group actions on  pseudo Riemannian hypersurfaces} \label{transactSec}
Now we explain how it is possible to use results from invariant
theory for reductive algebraic groups to classify pseudo Riemannian
hypersurfaces of degree $d$ admitting a transitive reductive algebraic group
of linear transformations. We will give a complete classification
for (quadratic and) cubic hypersurfaces, which is the case 
relevant to homogeneous
special geometry. 

Let $G$ be a real algebraic reductive group and $V$ an (algebraic) $G$-module. 
Assume that the connected component $G_0$ acts transitively on the hypersurface
${\cal H}\subset {\cal H}_1(h) \subset V$, where $h$ is the basic polynomial,
s.\ Def.\ \ref{basicDef}. Then $G$ preserves $h$, i.e.\ $h$ is a 
$G$-invariant. In fact, any $G$-invariant is of the form $ch^k$, $c\in 
\mbox{\Bbb C}$, $k \in \mbox{\Bbb N}_0$. 

Consider now the one dimensional
extension $G^1 = \mbox{\Bbb R}^* \times G$ and on $V$ the canonical 
structure of $G^1$-module, where  $\mbox{\Bbb R}^*$ acts by the 
standard scalar multiplication on $V$.  The $G^1$-module $V$ has an open
orbit, it is a {\bf prehomogeneous vector space} (P.V.) in the terminology
of \cite{S-K}. 

In \cite{S-K} T.\ Kimura and M.\ Sato have classified all {\em 
irreducible} P.V.s for {\em complex} algebraic reductive groups. 
If we consider the $(G^1)^{\mbox{\Bbb C}}$-module $V^{\mbox{\Bbb C}}$ 
and assume that it is irreducible, then it must appear in the classification  
\cite{S-K}. The $(G^1)^{\mbox{\Bbb C}}$-module $V^{\mbox{\Bbb C}}$  is 
irreducible if and only if the $G^{\mbox{\Bbb C}}$-module $V^{\mbox{\Bbb C}}$
is irreducible. We will assume first that this condition is satisfied and treat
the reducible case later. Without restriction of generality we also assume 
that $G^{\mbox{\Bbb C}}$ acts almost faithfully, i.e.\ with only discrete 
kernel. Then $V^{\mbox{\Bbb C}}$ is an irreducible almost faithful P.V.\ of 
$(G^1)^{\mbox{\Bbb C}}$ and by a theorem of Cartan (Thm.\ 1 of \cite{S-K}) 
$G$ must be 
semisimple. Remark that since $h$ is a $G$-invariant, it must be a {\bf 
relative} $G^1$- and $(G^1)^{\mbox{\Bbb C}}$-{\bf invariant}, i.e.\ 
$h$ is preserved up to scaling. A P.V. was called {\bf regular} in \cite{S-K}
if there exists a relative invariant with not identically vanishing
Hessian determinant. Remark that a one dimensional P.V. is always 
regular. 
\begin{lemma}
If  ${\cal H}= G_0v\subset {\cal H}_1(h) \subset V$ is a pseudo Riemanian
hypersurface in the $G$-module $V\ni v$, then $V^{\mbox{\Bbb C}}$  is 
a regular P.V. of $(G^1)^{\mbox{\Bbb C}}$.  
\end{lemma}  

\noindent
{\bf Proof:} If $\partial^2h|T_v{\cal H}$ is nondegenerate, then 
${\partial^2 h|}_v$ is  nondegenerate, since $h$ is a homogeneous 
polynomial of degree $d\ge 2$. $\Box$ 

An analogeous discussion applies to complex Riemannian hypersurfaces 
${\cal H} \subset V \cong \mbox{\Bbb C}^n$ admitting a transitive irreducible
linear action of a complex algebraic reductive group $G$. We can 
easily deduce the classification of such group actions from 
the classification of P.V.s in \cite{S-K}. We will consider two 
representations $R :  G \rightarrow GL(V)$ and 
$R' :  G' \rightarrow GL(V')$ as {\bf equivalent} if there is an isomorphism
$GL(V) \stackrel{\sim}{\rightarrow} GL(V')$ mapping $R(G)$ onto $R(G')$. 

\begin{thm} \label{irredcxredmodThm} The following list gives, up to 
equivalence, all irreducible 
$G$-modules $V$ of a 
connected complex algebraic reductive group $G$ which induce a transitive 
action on a complex Riemannian cubic hypersurface ${\cal H} \subset 
{\cal H}_1(h)\subset V$. The quadruples $(V^n,G,h,K)$ 
below contain the $G$-module $V$ of (complex) dimension $n$, the basic  
cubic polynomial $h$ (unique up to multiplicative constant) and the 
isotropy group $K$ of a point $v\in {\cal H}$ as abstract group.\\  
1) $(V^9 = U\otimes \mbox{\Bbb C}^3,\; H \times SL(3,\mbox{\Bbb C}),\; \det ,\;
 H)$,
where $U = \mbox{\Bbb C}^3$ is the standard 3-dimensional module of 
$H = SL(3,\mbox{\Bbb C})$, 
$H = SO(3,\mbox{\Bbb C})$ or $H = \{ e\}$ and 
$\mbox{\Bbb C}^3$ is the standard $SL(3,\mbox{\Bbb C})$-module. \\
2) $(V^6 = \vee^2 \mbox{\Bbb C}^3,\; SL(3,\mbox{\Bbb C}),\; \det ,\;
SO(3,\mbox{\Bbb C}))$,\\
3) $(V^{15} = \wedge^2 \mbox{\Bbb C}^6,\; SL(6,\mbox{\Bbb C}),\; Pff, \; 
Sp(3,\mbox{\Bbb C}))$, where $Pff$ is the Pfaffian of a skew symmetric
$6\times 6$-matrix: 
\[ Pff (A) = \sum_{\sigma \in S_6} {sgn}(\sigma ) a_{\sigma (1) \sigma (2)} 
 a_{\sigma (3) \sigma (4)}  a_{\sigma (5) \sigma (6)}\, , \quad 
A = (a_{ij})\, , \; i,j= 1,\ldots ,6\, .\] 
(Here we use the notation $Sp(n,\mbox{\Bbb C})$ for the symplectic group of 
$\mbox{\Bbb C}^{2n}$ as in \cite{S-K}.)      \\
4)  $(V^{27} = {Herm}_3(\mbox{\Bbb O}) \otimes \mbox{\Bbb C} = 
{Herm}_3({\cal C})\otimes \mbox{\Bbb C} ,\; E_6,\; \det 
,\; F_4)$, where ${Herm}_3(\mbox{\Bbb O})$ (resp.\ ${Herm}_3({\cal C})$)
denotes the Hermitian 
$3\times 3$ matrices over the octonions (resp.\ over the real split Cayley
algebra $\cal C$), which is an irreducible module of the
real form $E_6^{(-26)}$ (resp.\ $E_6^{(6)}$) of $E_6$. 
\end{thm} 

Now we classify all real $G$-modules $V$ of real algebraic reductive groups
$G$ such that the connected component $G_0$ acts transitively on a cubic 
pseudo Riemannian hypersurface ${\cal H} \subset {\cal H}_1(h) \subset V$ 
under the condition that the $G^{\mbox{\Bbb C}}$-module $V^{\mbox{\Bbb C}}$
is irreducible. If the last condition is satisfied, we shall say that
the module $(V,G)$ is {\bf totally irreducible}. It is clear that under
these conditions the complexified module $(V^{\mbox{\Bbb C}},
G^{\mbox{\Bbb C}})$ must be one of the list in Theorem \ref{irredcxredmodThm}. 
In other words, we must study the real forms of the complex modules 
classified in Theorem \ref{irredcxredmodThm}. 

Let now $V$ be a complex $G$-module of a connected complex algebraic reductive
group $G$. 

\begin{dof} A {\bf real structure} for the module $(V,G)$ is given by 
\begin{itemize}
\item[1)] a (real algebraic) antiholomorphic automorphism $\tau : G
\rightarrow G$ of the complex algebraic group $G$ and 
\item[2)] a complex antilinear involution $\tau : V\rightarrow V$ of the 
complex vector space $V$ (denoted by the same letter) such that
\[ \tau (av) = \tau (a) \tau (v) \, , \quad a\in G\, , \quad v\in V\, .\]  
\end{itemize} 
The pair $(V^{\tau}, G^{\tau})$ is called a {\bf real form} for $(V,G)$, 
where the superscript $\tau$ stands for fixed point set of $\tau$.   
\end{dof} 

Now we assume that ${\cal H} = Gv \subset {\cal H}_1(h)$ is a complex 
Riemannian hypersurface for some $v\in V$. Remark that then $Gv$
is a complex Riemannian hypersurface for all $v\in V-S$, where the singular 
set $S$ is Zariski closed. Given a real structure $\tau$ for $(V,G)$ we may 
assume
that $h$ is real, i.e.\ $\overline{h(v)} = h(\tau v)$ for all $v\in V$. 
In fact, multiplying $h$ by a suitable complex constant we may assume
that $h^{\tau}: = {Re}\, {h|}_{V^{\tau}} \neq 0$. Then the complex 
extension $h^{\tau}_{\mbox{\Bbb C}}$ is a $G$-invariant polynomial 
on $V$ of the same degree as $h$ and hence $h^{\tau}_{\mbox{\Bbb C}} = 
ch$ for some $c\in \mbox{\Bbb C}^{\ast}$. Finally,  
we can assume that ${\cal H} = Gv$ for some $v\in V^{\tau}$. Indeed,   
the singular set $S$ is a Zariski closed proper subset of $V$ and hence 
$V^{\tau}\not\subset S$. Under these assumptions we have the following lemma. 

\begin{lemma} The real structure $\tau$ preserves the isotropy group
$K$ of $G$ at $v$ and the real form $K^{\tau}$ of $K$ is precisely 
the isotropy group of $v$ in $G^{\tau}$. Moreover, the orbit 
$(G^{\tau})_0v \cong (G^{\tau})_0/(K^{\tau}\cap (G^{\tau})_0)$ of the 
connected component $(G^{\tau})_0$ of $G^{\tau}$ coincides 
with the connected component $({\cal H}^{\tau})_0$ of $v$ in the fixed
point set ${\cal H}^{\tau}$, ${\cal H} = Gv$. Finally,  $({\cal H}^{\tau})_0$
is pseudo Riemannian if and only if ${\cal H}$ is complex Riemannian. 
\end{lemma} 

\noindent
{\bf Proof:} It is clear that  $K^{\tau} = K\cap G^{\tau}$ is the isotropy group of $v$ in 
$G^{\tau}$; in particular $G^{\tau}v \cong G^{\tau}/K^{\tau}$ and 
$(G^{\tau})_0v \cong (G^{\tau})_0/(K^{\tau}\cap (G^{\tau})_0)$. From 
$v\in V^{\tau}$ it follows that $\tau K = K$, $\tau (Gv) = Gv$ and 
$G^{\tau}v \subset (Gv)^{\tau}$, so $K^{\tau}$ is a real form of 
$K$ and 
\[ \dim_{\mbox{\Bbb C}} Gv \ge \dim_{\mbox{\Bbb R}}(Gv)^{\tau} \ge 
\dim_{\mbox{\Bbb R}}G^{\tau}v \] 
\[ = \dim_{\mbox{\Bbb R}}G^{\tau}/K^{\tau} 
= \dim_{\mbox{\Bbb C}} G/K = \dim_{\mbox{\Bbb C}} Gv \, .\] 
This implies $\dim (Gv)^{\tau} = \dim G^{\tau}v$ and $(G^{\tau})_0v =
((Gv)^{\tau})_0$. The last statement of the lemma is trivial since we are 
assuming that the basic polynomial $h$ is real. $\Box$ 

\begin{thm} \label{irredrealredThm}The totally irreducible modules of real 
algebraic reductive 
groups inducing a transitive action of the group's identity component on a 
cubic pseudo Riemannian hypersurface
are (up to equivalence) the  real forms $(V^{\tau}, G^{\tau})$ 
of the complex modules $(V,G)$ classified in Theorem \ref{irredcxredmodThm}.  
We list the triples $(V^{\tau}, G^{\tau}, K^{\tau})$, where $K^{\tau}$ is 
the isotropy group of $v\in V^{\tau}$; ${\cal H} = Gv\subset V$  the 
underlying complex Riemannian hypersurface. The corresponding homogeneous
pseudo Riemannian cubic hypersurfaces are all locally symmetric.
(The real cubic hypersurfaces with positively or negatively defined canonical
metric  correspond to the triples  $(V^{\tau}, G^{\tau}, K^{\tau})$ 
with compact isotropy group $K^{\tau}$.) \\
1)$(U^0\otimes \mbox{\Bbb R}^3,\; H^0 \times SL(3,\mbox{\Bbb R}),\; 
 H^0)$,
where $U^0 = \mbox{\Bbb R}^3$ is the standard 3-dimensional module of 
$H^0 = SL(3,\mbox{\Bbb R})$, 
$SO(1,2)$, $SO(3)$   or 
$\{ e\}$. \\ 
2)$(\vee^2 \mbox{\Bbb R}^3,\; SL(3,\mbox{\Bbb R}),\; SO(3))$,\\ 
3) $( \wedge^2 \mbox{\Bbb R}^6,\; SL(6,\mbox{\Bbb R}),\; Sp(3,\mbox{\Bbb R}))$,
$((\wedge^2_{\mbox{\Bbb C}}\mbox{\Bbb H}^3)^{\tau},\; SL(3,\mbox{\Bbb H}),\; 
Sp(3))$, where the real structure $\tau$ on 
$\wedge^2_{\mbox{\Bbb C}}\mbox{\Bbb H}^3 \subset \mbox{\Bbb H}^3
\otimes_{\mbox{\Bbb C}} \mbox{\Bbb H}^3$ is the square of the 
quaternionic structure on the complex vector space $\mbox{\Bbb H}^3$.\\ 
4)  $({Herm}_3(\mbox{\Bbb O}),\; E_6^{(-26)},\; F_4^{(-52)})$, 
$({Herm}_3({\cal C}),\;$ $E_6^{(6)},\;$ 
$F_4^{(4)})$, for the notation s.\ Theorem \ref{irredcxredmodThm}, 4). 
\end{thm} 

\noindent
{\bf Proof:} The fact that the above list gives all real forms 
for the complex modules of Theorem \ref{irredcxredmodThm} can be 
checked using Tits' tables \cite{Tt}. $\Box$ 

Now we study {\em reducible} $G$-modules of connected complex algebraic
reductive groups $G$ inducing a transitive $G$-action on a nondegenerate
hypersurface. We  use the following theorem as lemma.  
\begin{lemma}\cite{S-K}\label{regL} Let $(V,G)$ be a prehomogeous vector space 
of a connected complex algebraic reductive group $G$. Then the following
conditions are equivalent.
\begin{enumerate}
\item[(i)] $(V,G)$ is regular.
\item[(ii)] The generic isotropy group $G_v$, $v\in V$, is reductive. 
\item[(iii)] The singular set $S = V- Gv$ is a hypersurface.  
\end{enumerate}
\end{lemma}

\begin{thm} \label{3casesThm} Let $(V,G)$ be a faithful module of a 
simply connected 
and complex algebraic reductive group $G$ acting transitively
on a nondegenerate cubic hypersurface ${\cal H} \subset V$ with basic
polynomial $h$. Consider its canonical extension $(V, G^1)$ to a   
regular prehomogeneous vector space $V$ of the group $G^1 = 
\mbox{\Bbb C}^{\ast}\times G$. Then $(V, G^1)$ is the sum of at most 
3 irreducible regular P.V.s. Moreover, only the three following 
possibilities can occur.
\begin{itemize}
\item[(1)] The module $(V, G^1)$ and the basic polynomial $h$ are irreducible
and $G$ is semisimple. 
\item[(2)] The $G^1$-module $V$ is the direct sum $V = V_1 \oplus V_2$ 
of two irreducible P.V.s and $\dim V_1 = 1$. The polynomial $h$ is the 
product $h = lq$ of a linear function $l$ on $V_1$ and a quadratic 
$G$-invariant $q$ on $V_2$. Finally, $G \cong \mbox{\Bbb C}^{\ast}\times G'$,
where $G'$ is semisimple or trivial. 
\item[(3)] The  $G^1$-module $V$ is the direct sum of three one-dimensional
P.V.s and $G\cong \mbox{\Bbb C}^{\ast}\times  \mbox{\Bbb C}^{\ast}$. 
\end{itemize} 
\end{thm}

\noindent
{\bf Proof:}  First we remark that the  $G^1$-module $V$ is the direct sum
$V = \oplus_{j=1}^r V_j$ of irreducible submodules, because $G^1$ is 
a reductive algebraic group. Moreover, the irreducible summands $V_j$ are
again prehomogeneous vector spaces. 

From the assumption that $G$ acts transitively on a {\em nondegenerate} 
hypersurface it follows that $(V,G^1)$ is a {\em regular} P.V. and 
by Lemma \ref{regL} the generic isotropy group $G_v$ is reductive.
This implies that the generic isotropy group for the irreducible summands 
$V_j$ is also reductive and hence the $V_j$ are regular P.V.s. In particular, 
each of them admits a non-constant relative invariant $h_j$. 

We have already remarked that  
$(V, G^1)$ being irreducible and faithful implies the semisimplicity of $G$.  
Now we prove that if $(V, G^1)$ is irreducible, then the basic polynomial
$h$ is also irreducible. In fact, assume that $h = fg$ is the product
of two non-constant polynomials. Then $f$ and $g$ are relative $G^1$-invariants
and one of them must be linear. However, a module $V$ with a linear relative
invariant $l$ splits as $V = V_1 \oplus V_2$, where $\dim V_1 =1$ and 
$V_2 = {ker} \, l$. 

Any  simply connected complex algebraic reductive group $G$
is the direct product $G = Z \times G'$ of its centre $Z$, which is an 
algebraic torus, and a semisimple group $G'$. 

If the $G^1$-module $V
= \oplus_{j=1}^r V_j$  is the sum of $r$ irreducible submodules $V_j$, 
then it has $r$ algebraically independent relative $G^1$-invariants
$h_1,\ldots ,h_r$, which are necessarily $G'$-invariant. Remark that at a
generic point $v\in V$ the differentials $dh_1,\ldots , dh_r$ are 
linearly independent. 
Since $G$ acts transitively on a hypersurface in $V$, the group
$G'$ has an orbit of codimension $\le \dim Z + 1$. This implies 
$r\le \dim Z + 1$. On the other hand, since $\mbox{\Bbb C}^{\ast} \times Z$ 
acts faithfully on $V$, we have that $1+ \dim Z \le r$, hence
\[ r \quad = \quad \dim Z + 1\, .\] 
Now we show that $r\le 3$ and hence $\dim Z = r-1 \le 2$. A homogeneous 
polynomial $f = f_1 \cdots f_r$ on $V$ which is a product of homogeneous
polynomials $f_j$ of degree $d_j$ on $V_j$ will be called a 
{\bf monomial of degree} $(d_1,\ldots , d_r)$, with respect to the 
$G$-invariant decomposition  $V= \oplus_{j=1}^r V_j$.  Any homogeneous 
polynomial $f$ on $V$ of degree $d$ can be decomposed into monomials
$f_{\delta}$ of degree $\delta = (d_1,\ldots , d_r)\in \mbox{\Bbb N}^r_0$, 
${|\delta|} = d_1 + \cdots + d_r = d$:
\[ f\quad = \quad \sum_{|\delta | =d} f_{\delta} \, .\] 
If $f$ is $G$-invariant, then  the $f_{\delta}$  are $G$-invariant. 
It follows that the basic cubic polynomial $h$ is a monomial 
$h = h_1 \cdots h_r$ and the degrees $d_j\ge 0$ of the homogeneous polynomials 
$h_j$ must add up to three: $d_1 + \cdots + d_r =3$. Now we use the 
fact that $h$ has nondegenerate Hessian form $\partial^2 h$ to conclude
that $d_j\ge 1$ for all $j = 1,\ldots ,r$. This implies $r\le 3$, as claimed
above.  

If $r=1$, $V$ and hence $h$ is irreducible. If $r=2$, $h = h_1h_2$, 
where we can assume that $h_1 = l$ is a linear polynomial on $V_1$ and 
$h_2 =q$ is a quadratic polynomial on $V_2$. From the irreducibility of 
$V_2$ (and nondegeneracy of $\partial^2 h$) it follows that either $q$
is irreducible (and a nondegenerate quadratic form) or it is the square
of a linear polynomial and $\dim V_2 =1$. The second case implies 
$G' = \{ e\}$ by the faithfullness of the $G$-action. 
Finally if $r=3$, $h = h_1 h_2 h_3$ is the product of three relative
$G^1$-invariant polynomials $h_j$ on $V_j$. In particular, 
$\dim V_j=1$ and $G' = \{ e\}$. $\Box$. 

Theorem \ref{3casesThm} reduces the classification of homogeneous pseudo 
Riemannian cubic hypersurfaces of reductive real algebraic groups $G$ to three 
cases. The first case is that of totally irreducible $G$-modules and was 
treated in Theo\-rem \ref{irredrealredThm}. The third case corresponds to 
the real forms of the flat complex Riemannian cubic hypersurface
in $\mbox{\Bbb C}^3$ with basic polynomial $h(z^1,z^2,z^3) = z^1z^2z^3$ and 
transitive action of the algebraic torus $\mbox{\Bbb C}^{\ast} \times 
\mbox{\Bbb C}^{\ast}$. The second case reduces to the classification of 
homogeneous pseudo Riemannian quadratic hypersurfaces, which we give now.

Every homogeneous quadratic polynomial $q$ on $\mbox{\Bbb K}^n$ such that
$\partial^2 q$ is nondegenerate is of the form 
\[ q(X) = {\langle X,X\rangle}\, , \quad X\in \mbox{\Bbb K}^n\, ,\] 
for some nondegenerate symmetric {\Bbb K}-bilinear form $\langle \cdot , 
\cdot \rangle$ on 
$\mbox{\Bbb K}^n$. If $\mbox{\Bbb K}=\mbox{\Bbb C}$, we call ${\cal H} =
{\cal H}_1(q)$ a {\bf complex sphere}. If $\mbox{\Bbb K}=\mbox{\Bbb R}$, 
the connected components ${\cal H}$ of ${\cal H}_1(q)$ are called 
{\bf pseudo spheres}. The canonical metric $g$ of $\cal H$ at $X_0 \in 
{\cal H}$ is precisely
\[ g_{X_0}(X,X) = - {\langle X,X\rangle}\, , \quad X\in T_{X_0}{\cal H} = \{ X\in 
 \mbox{\Bbb K}^n| {\langle X_0,X\rangle } = 0\}\, .\] 
In the real case this shows that $({\cal H},g)$ has signature $(l,k-1)$ if 
$\langle \cdot , \cdot \rangle$ has signature $(k,l)$. In particular, the canonical 
metric of the hyperboloid ${\cal H}_1(q) \subset \mbox{\Bbb R}^{1,n-1}$ 
is positively defined and that of the sphere  ${\cal H}_1(q) \subset 
\mbox{\Bbb R}^{n,0}$ is negatively defined; where $q(X) = x_1^2 + \cdots 
x_k^2 - x_{k+1}^2 - \cdots -x_n^2$ for $\mbox{\Bbb R}^{k,l}$, $k+l = n$. 

By Witt's Theorem, $SO_0(k,l)$ acts transitively on the (connected) 
pseudo spheres ${\cal H} \subset \mbox{\Bbb R}^{k,l}$. The identity 
component of the isotropy group of $SO_0(k,l)$ at $(1,0, \ldots , 0)$ 
is $SO_0(k-1,l)$. Moreover, the canonical metric  defines
on $\cal H$ the structure of pseudo Riemannian symmetric space. 
The symmetry $\sigma_{X_0}: {\cal H} \rightarrow {\cal H}$ at 
$X_0 \in {\cal H}$ is given by: 
\[  \sigma_{X_0}(X) = -\left( X - \frac{\langle X,X_0\rangle}{\langle X_0,X_0
\rangle}X_0\right) + 
\frac{\langle X,X_0\rangle}{\langle X_0,X_0\rangle}X_0\, .\] 
We sum up our discussion. 

\begin{prop} Any pseudo Riemannian quadratic hypersurface (always
with the canonical metric) is a pseudo sphere ${\cal H}\subset 
\mbox{\Bbb R}^{k,l}$   and, in particular, a pseudo Riemannian 
symmetric space. 
Its full connected group of linear automorphisms is $SO_0(k,l)$, which acts 
irreducibly on  $\mbox{\Bbb R}^{k,l}$. 
\end{prop} 

To round up our discussion, we give the classification of transitive
linear algebraic actions of reductive groups on complex and on pseudo 
Riemannian quadratic hypersurfaces. As before, the complex case can be easily 
extracted from \cite{S-K}:
\begin{thm} \label{redcxThm} Any $G$-module $V$ of a complex algebraic 
reductive group
which induces a transitive action on a nondegenerate quadratic 
hypersurface is (up to equivalence) one in the following list. All of them
are irreducible and hence their basic quadratic polynomial $q$ is 
the (up to scaling) unique $G$-invariant quadratic form on $V$. 
The corresponding 
quadratic hypersurfaces are complex spheres. The last entry $K$ in the 
triples $(V^n,G,K)$ is, as before, the generic isotropy group as abstract
group and $n = \dim V$. \\
1) $(V^4 = U \times \mbox{\Bbb C}^2,\; H \times SL(2,\mbox{\Bbb C}),\; H)$, 
$U = \mbox{\Bbb C}^2$, $H = SL(2,\mbox{\Bbb C})$ or $=\{ e\}$, \\
2) $(V^4 = \vee^2 \mbox{\Bbb C}^2,\; SL(2,\mbox{\Bbb C}),\; 
SO(2,\mbox{\Bbb C}))$, \\
3) $(V^{4n} =  \mbox{\Bbb C}^{2n} \otimes  \mbox{\Bbb C}^2,\; 
Sp(n, \mbox{\Bbb C}) \times SL(2, \mbox{\Bbb C}),\; 
Sp(n-1, \mbox{\Bbb C}) \times Sp(1, \mbox{\Bbb C}))$, $n\ge 2$. 
In this case the basic polynomial $q$ can be described as follows. Consider 
$\mbox{\Bbb C}^{2n} \otimes  \mbox{\Bbb C}^2$ as vector space 
of  complex $(2n) \times 2$-matrices $A = (a_{ij})$, then 
\[ q(A) = Pff (A^t J A)\, , \quad J = 
\left(\begin{array}{rr} 0 & \mbox{\boldmath$1_n$}\\ 
  \mbox{\boldmath$-1_n$}  & 0
\end{array} \right)\, ,\]
  is the Pfaffian of the skew
symmetric matrix $A^t J A$. \\
4) $(V^n = \mbox{\Bbb C}^n,\; SO(n, \mbox{\Bbb C}),\; 
SO(n-1,\mbox{\Bbb C}))$,\\ 
5) $(V^8 = \mbox{spinor module},\; Spin(7,\mbox{\Bbb C}),\; G_2)$, \\
6) $(V^{16} = \mbox{spinor module},\; Spin(9,\mbox{\Bbb C}),\;
Spin(7,\mbox{\Bbb C}))$, \\ 
7) $(V^7 \mbox{with highest weight $\Lambda_2$},\; G_2,\; 
SL(3,\mbox{\Bbb C}))$. 
\end{thm}

\begin{thm} The real algebraic $G$-modules $V$ of reductive groups $G$ 
with transitive action of $G_0$ on a quadratic pseudo Riemannian hypersurface 
${\cal H} \subset V$ are obtained from the real forms of the complex 
modules in Theorem \ref{redcxThm} and are listed below.    
The hypersurface  ${\cal H} \subset {\cal H}_1(q) \subset (V,q)$ is a 
pseudo sphere in the 
pseudo  Euclidean vector space $(V,q)$, where $q$ is a real 
basic $G$-invariant quadratic polynomial of the irreducible $G$-module
$V$.\\  
1) $(\mbox{\Bbb R}^2 \times \mbox{\Bbb R}^2,\; 
H^0 \times SL(2,\mbox{\Bbb R}))$, where $H^0 = SL(2,\mbox{\Bbb R})$ or 
$=\{ e\}$, $(\mbox{\Bbb R}^4,\; SO(4))$, cf.\ 4), \\
2) $(\vee^2 \mbox{\Bbb R}^2,\; SL(2,\mbox{\Bbb R}))$, \\
3) $(\mbox{\Bbb R}^{2n} \otimes  \mbox{\Bbb R}^2,\; 
Sp(n, \mbox{\Bbb R}) \times SL(2, \mbox{\Bbb R}))$, 
$(\mbox{\Bbb H}^n,\; Sp(n-l,l) \times Sp(1))$, $0\le l\le n\ge 2$, \\  
4) $(\mbox{\Bbb R}^{k,l},\; SO(k,l ))$, $k+l =n$,\\ 
5) spinor module of $Spin(7,0)$ and $Spin(3,4)$ and semi spinor modules
of $Spin(0,7)$ and $Spin(4,3)$,\\ 
6) spinor module of $Spin(8,1)$, $Spin(4,5)$ and $Spin(0,9)$ and 
semi spinor modules of $Spin(1,8)$, $Spin(5,4)$ and $Spin(9,0)$, \\ 
7) $({Im}\, \mbox{\Bbb O},\; G_2^{(-14)} = {Aut}\, \mbox{\Bbb O})$,  
$({Im}\, {\cal C},\;$ $G_2^{(2)} = {Aut}\, {\cal C})$, 
where ${Im}\, \mbox{\Bbb O}$ (resp.\ 
${Im}\, {\cal C}$) denotes the imaginary part and ${Aut}\, \mbox{\Bbb O}$ 
(resp.\ ${Aut}\, {\cal C}$) the full irreducible automorphism group of the 
octonions (resp.\ of the real split Cayley algebra).   
\end{thm} 

\noindent
{\bf Proof:} The real forms are obtained again using tables, s.\ 
\cite{Tt} and \cite{L-M}. $\Box$ 
\subsection{Appendix: Special K\"ahler submanifolds of Alekseevsky's 
quaternionic K\"ahler manifolds and their pseudo Riemannian analogues} 
\label{qKSec}
First of all we explain how one can obtain pseudo K\"ahlerian versions of the 
cubic (and also of the quadratic) normal J-algebras classified in Theorem 
\ref{classhomogHThm}.
Consider e.g.\ the cubic normal J-algebras $({\fr u}_0 (p),\langle \cdot ,
\cdot \rangle ,J)$ of rank 2. If $p\neq 0$, we can define a new  pseudo 
K\"ahler Lie algebra $({\fr u}_0 (p)',\langle \cdot ,
\cdot \rangle ',J)$ changing the scalar product on the subspace 
${\fr x}_1 \subset {\fr u}_0 (p)$ only by a sign. The formulas defining
the Lie bracket on ${\fr u}_0 (p)'$ are the same as for ${\fr u}_0 (p)$
only the scalar product occuring in these formulas is substituted by
$\langle \cdot ,\cdot \rangle '$. 

Recall (Theorem \ref{PSThm}) that to every isometric map 
$\psi : {\fr x}_{23}^- \times {\fr x}_{12}^- 
\rightarrow {\fr x}_{13}^-$ of Euclidean vector spaces ${\fr x}_{23}^-$, 
${\fr x}_{12}^-$ and ${\fr x}_{13}^-$ we can associate a normal  
J-algebra ${\fr u}_0 (\psi )$ of type I. If moreover ${\fr x}_{23}^- = 0$ or if 
$\dim {\fr x}_{12}^- = \dim {\fr x}_{13}^-$, then ${\fr u}_0 (\psi )$ is a
cubic normal J-algebra, s.\ Theorem \ref{classhomogHThm}. Similarly, 
to every isometric map 
$\psi : {\fr x}_{23}^- \times {\fr x}_{12}^- 
\rightarrow {\fr x}_{13}^-$  of 
{\em pseudo}    Euclidean vector spaces ${\fr x}_{23}^-$, 
${\fr x}_{12}^-$ and ${\fr x}_{13}^-$
we can associate a {\em pseudo} K\"ahlerian Lie  algebra ${\fr u}_0 (\psi )$
and we have the following proposition. 

\begin{prop} \label{pssKProp} The pseudo K\"ahler Lie algebra 
$({\fr u}_0 (p)',\langle \cdot ,
\cdot \rangle ',J)$ is cubic.  Let  $\psi : {\fr x}_{23}^- 
\times {\fr x}_{12}^- 
\rightarrow {\fr x}_{13}^-$ be an isometric map  of 
pseudo   Euclidean vector spaces ${\fr x}_{23}^-$, 
${\fr x}_{12}^-$ and ${\fr x}_{13}^-$ and assume that 
$\psi$ is  special or of order $0$, s.\ Definition \ref{sim}. Then the pseudo 
K\"ahlerian Lie algebra $({\fr u}_0 (\psi ),
\langle \cdot ,\cdot \rangle ,J)$ associated to $\psi$ is cubic. 
In particular,
$({\fr u}_0 (p)',\langle \cdot ,\cdot \rangle ',J)$ and $({\fr u}_0 (\psi ),
\langle \cdot ,\cdot \rangle ,J)$ are 
pseudo K\"ahler Lie algebras for   special pseudo
K\"ahler tube domains. 
\end{prop}

\noindent
{\bf Proof:} The proof is analogeous to the proof of Theorem 
\ref{classhomogHThm}. 
$\Box$ 

\medskip
Next we recall some basic facts about Alekssevsky's quaternionic
K\"ahler manifolds, s.\ \cite{A}, \cite{C} and \cite{A-C} for 
details. 
\begin{dof} An {\bf Alekseevsky space} is a quaternionic K\"ahler manifold
$M$ admitting a simply transitive (non-Abelian) splittable solvable 
group $\cal L$ of isometries.
\end{dof} 
We can present $M$ as metric Lie group $({\cal L},g)$ and consider
its metric Lie algebra $({\fr l}, \langle \cdot , \cdot \rangle )$.
The quaternionic K\"ahler structure of $M$ induces a quaternionic
structure ${\fr q} = {span} \{ J_1,J_2,J_3\}$ on the Euclidean 
vector space $({\fr l}, \langle \cdot , \cdot \rangle )$. The triple
$({\fr l}, \langle \cdot , \cdot \rangle ,{\fr q})$ associated
to the Alekseevsky space $M$ is called its {\bf Alekseevskian
Lie algebra}. According to \cite{A}, \cite{dW-VP} and \cite{C} 
there are (up to symmetric spaces) 3 series of Alekseevsky spaces:
$\cal T$, $\cal W$ and $\cal V$-spaces. Their Alekseevskian Lie
algebras $({\fr l}, \langle \cdot , \cdot \rangle ,{\fr q})$
are constructed as quaternionic K\"ahler extensions of normal J-algebras
$\fr u$. 

More precisely, given a normal J-algebra of the form 
${\fr u} = {\fr f}_0 \stackrel{\perp}{\oplus} {\fr u}_0$ which 
admits a so called {\bf Q-representation} $T: {\fr u} \rightarrow
{End} (\tilde{\fr u})$ we can canonically define the structure
$({\fr l}, \langle \cdot , \cdot \rangle ,{\fr q})$ of 
Alekseevskian Lie algebra   on the vector space 
${\fr l} = {\fr u} + \tilde{\fr u}$ such that 
$({\fr u},\langle \cdot , \cdot \rangle |{\fr u}, J=J_1|{\fr u})$ is a K\"ahlerian subalgebra, $\tilde{\fr u} = J_2{\fr u}$, $[\tilde{\fr u},
\tilde{\fr u}] \subset {\fr u}$ and $[ {\fr u} , \tilde{\fr u}] \subset
\tilde{\fr u}$ is given by $T$. We remark that the definition of 
Q-representation and this construction can be naturally generalized to 
pseudo K\"ahlerian Lie 
algebras; the homogeneous spaces associated to  
$({\fr l}, \langle \cdot , \cdot \rangle ,{\fr q})$ being {\em pseudo} 
quaternionic K\"ahler manifolds. 

The 3 series of  $\cal T$, $\cal W$ and $\cal V$-spaces are defined 
by 3 series of normal J-algebras $({\fr u}_0,\langle \cdot , \cdot \rangle ,J)$
and corresponding Q-representations. 
By direct comparison of these series of normal J-algebras ${\fr u}_0$ with the 
cubic normal J-algebras (s.\ Theorem \ref{classhomogHThm}) one can easily check
the following facts. 

The Alekseevsky spaces ${\cal V}(\psi)$ are defined by the 
Q-representation
of the  normal J-algebras ${\fr u} = {\fr f}_0 \oplus  
{\fr u}_0(\psi )$ for which ${\fr u}_0(\psi )$ is a cubic normal
J-algebra of rank 3 defined by a non zero special isometric map $\psi$.

For the spaces ${\cal T}(p)$, $p = 0,1,2,\ldots$, 
\[ {\fr u}_0(p) = {\fr e}(p+1,1) {\in\!\!\!\!\!+} 
{\fr e}(1,\frac{1}{\sqrt{2}}) \] 
are the cubic normal J-algebras
of rank 2, cf.\ Proposition \ref{rk=2deg=dProp}. 

Finally, for  the spaces 
${\cal W}(p,q) \cong {\cal W}(q,p)$, $p,q=0,1,2,\ldots$, 
\[ {\fr u}_0(p,q) = {\fr e}(p+q+1,1) {\in\!\!\!\!\!+} 
{\fr e}(1,1){\in\!\!\!\!\!+} {\fr e}(1,1) \]
are the cubic normal J-algebras of
rank 3 associated to isometric maps of order zero. 

The preceding facts motivate the following definition.
\begin{dof} Let  $M = ({\cal L},g)$ be an Alekseevsky space with 
Alekseevskian Lie algebra $({\fr l}, \langle \cdot , \cdot \rangle ,{\fr q})$, 
${\fr l} = {\fr u} + \tilde{\fr u}$, ${\fr u} =  
{\fr f}_0 \oplus {\fr u}_0$ the decompositions introduced
above and $({\cal U}_0,g,\tilde{J})$ the K\"ahler Lie group associated
to the normal J-algebra $({\fr u}_0,\langle \cdot , \cdot \rangle ,
J=J_1|{\fr u}_0)$. It is naturally identified with a (totally geodesic) 
K\"ahler submanifold
of $M$, which is called the {\bf special K\"ahler submanifold} of $M$.
The normal J-algebra $({\fr u}_0,\langle \cdot , \cdot \rangle ,
J=J_1|{\fr u}_0)$ is called the {\bf special K\"ahler subalgebra}
of the Alekseevskian Lie algebra  
$({\fr l}, \langle \cdot , \cdot \rangle ,{\fr q})$. 
  
\end{dof} 

We remark that to the cubic pseudo K\"ahlerian Lie algebras
of Proposition \ref{pssKProp} one can associate pseudo Riemannian analogues 
of the Alekseevsky spaces. 

\begin{prop} 
The cubic pseudo K\"ahlerian Lie algebras $({\fr u}_0 (p)',\langle \cdot ,
\cdot \rangle ',J)$ and ${\fr u}_0 (\psi )$
associated to a an isometric map $\psi$ of 
pseudo    Euclidean vector spaces, which is special or of zero order, admit 
a Q-representation. 
In particular, to any such map $\psi$ we can associate
a homogeneous pseudo quaternionic  K\"ahler manifold. 
\end{prop} 

\noindent
{\bf Proof:} The proof does not depend on the signature of the scalar
products, s.\ \cite {A}, \cite{C}. $\Box$ 

\medskip\noindent 
{\bf Remark 6:} Being a totally geodesic K\"ahlerian submanifold of a 
symmetric space,  the special K\"ahler submanifold of a symmetric 
Alekseevsky spaces is Hermitian symmetric. It can be described as follows. 
Let $M = G/K$ be a symmetric Alekseevsky space, $G$ its maximal connected
isometry group and $K = Sp(1)\cdot H$ its isotropy subgroup; $G$ and $H$
are  semisimple. We denote by 
$({\fr l}, \langle \cdot , \cdot \rangle ,{\fr q})$ the Alekseevskian 
Lie algebra corresponding to $M$ and  by 
$({\fr u}_0,\langle \cdot , \cdot \rangle ,J)$ its special K\"ahler subalgebra.
Then $\fr l$ is isomorphic to the Iwasawa Lie algebra of $G$ and ${\fr u}_0$
is isomorphic to the Iwasawa Lie algebra of a non compact semisimple 
Lie group $H^0\subset G$. The special K\"ahler submanifold $U \subset M$ is 
the orbit of the point $K = eK \in M = G/K$ under the group $H^0$ and is a 
Hermitian symmetric space of  non compact type.  

Consider the twistor space $Z = G/ (U(1)\cdot H)$ of $M$. It carries 
a natural structure of homogeneous complex contact manifold, s.\ \cite{W}. 
The contact hyperplane $D_S \subset T_SZ$, $S\in Z$, carries an 
$H$-invariant complex symplectic structure $\omega$. The representation of 
$H$ on $D_S$ is irreducible, preserves this complex symplectic structure 
and the orbit of the highest root 
vector is the base of a Lagrangean cone  $\hat{\cal C}$. Its projectivization 
$P(\hat{\cal C})$ is isomorphic to the compact dual $U^{\ast}$ of 
the Hermitian symmetric space $U$. 

We can consider the action of the complex semisimple Lie group 
$H^{\mbox{\C}}$ on the complex symplectic vector space $D_S$ and on the 
compact Hermitian symmetric space $U^{\ast} \cong P(\hat{\cal C})$. 
Now we remark that the Lie group $H^0 = {Isom} (U)$ is a (non compact)
real form of $H^{\mbox{\C}}$, hence we can also consider the action
of $H^0$ on $D_S$ and on $U^{\ast} \cong P(\hat{\cal C})$.
In this way we can realize $U$ as open orbit 
$U = P({\cal C}) \subset P(\hat{\cal C})$ of $H^0$ on the compact Hermitian
symmetric space $U^{\ast} \cong P(\hat{\cal C})$. 

One can check that $D_S$ carries a $H^0$-invariant
real structure $\tau$, such that the canonical special K\"ahler metric of 
the projectivized cone $P({\cal C})$ defined by the data $(D_S, \omega 
, \tau )$ coincides with the (nonpositively curved) Hermitian symmetric
metric given by the inclusion $U\subset M$. 

It is known\footnote{The author has learned this fact from J.-M.\ Hwang.}
that the compact projectivized cone $P(\hat{\cal C}) \cong U^{\ast}$
is precisely the projectivization of the cone $\hat{\cal C}
\subset T_SZ^{\ast}$ of 
tangent directions to minimal rational curves through a point $S\in Z^{\ast}$
in the twistor space $Z^{\ast}$ of the Wolf space $M^{\ast} = G^{\ast}/K$.
$M^{\ast} = G^{\ast}/K$ is the compact symmetric quaternionic 
K\"ahler manifold which is dual to the symmetric Alekseevsky space
$M = G/K$.  

\end{document}